\documentclass[12pt]{article}
\usepackage{amsmath, amsthm, amsfonts, stmaryrd}

\def\FF{\mathbb{F}}
\def\NN{\mathbb{N}}
\def\QQ{\mathbb{Q}}
\def\ZZ{\mathbb{Z}}

\def\bv{\mathbf{v}}
\def\bw{\mathbf{w}}

\numberwithin{equation}{subsection}
\newtheorem{theorem}[equation]{Theorem}
\newtheorem{cor}[equation]{Corollary}
\newtheorem{prop}[equation]{Proposition}
\newtheorem{lemma}[equation]{Lemma}
\newtheorem{conj}[equation]{Conjecture}
\theoremstyle{definition}
\newtheorem{defn}[equation]{Definition}

\newtheorem{remark}[equation]{Remark}

\DeclareMathOperator{\Gal}{Gal}

\begin{document}

\title{Finite automata and algebraic extensions of function fields}
\author{Kiran S. Kedlaya \\
Department of Mathematics \\
Massachusetts Institute of Technology \\
Cambridge, MA 02139 \\
\texttt{kedlaya@math.mit.edu}
}
\date{May 3, 2005}

\maketitle

\begin{abstract}
We give an automata-theoretic description of the algebraic closure of
the rational function field $\FF_q(t)$
over a finite field $\FF_q$, generalizing a result of Christol.
The description takes place within the Hahn-Mal'cev-Neumann field
of ``generalized power series'' over $\FF_q$.
Our approach
includes a characterization of well-ordered sets of rational numbers
whose base $p$ expansions are generated by a finite automaton,
as well as some techniques for computing in the algebraic closure;
these include an adaptation to positive characteristic
of Newton's algorithm for finding local expansions
of plane curves. We also conjecture a generalization of our results
to several variables.
\end{abstract}

\section{Introduction}

\subsection{Christol's theorem, and its limits}

Let $\FF_q$ be a finite field of characteristic $p$, and let
$\FF_q(x)$ and $\FF_q((t))$ denote the fields of rational functions
and of formal (Laurent) power series, respectively, over $\FF_q$.
Christol \cite{christol} (see also \cite{ckmr}) proved that an element
$x = \sum_{i=0}^\infty x_i t^i$ of $\FF_q((t))$
is algebraic over $\FF_q(t)$ (that is, is the root of a monic polynomial
in one variable with coefficients in $\FF_q(t)$) if and only if for each
$c \in \FF_q$, the set of base $p$ expansions of the integers $i$ for
which $x_i = c$ is generated by a finite automaton.

However, this is not the end of the story, for there are monic polynomials over
$\FF_q(t)$ which do not have any roots in $\FF_q((t))$, even if you enlarge
the finite field and/or replace $t$ by a root. An example, due to
Chevalley \cite{che}, is the polynomial
\[
x^p - x - t^{-1}
\]
Note that this is a phenomenon restricted to positive characteristic
(and caused by wild ramification): an old
theorem of Puiseux \cite[Proposition~II.8]{ser}
implies that if $K$ is a field of
characteristic 0, then any monic polynomial
of degree $n$ over $K(t)$ factors into linear polynomials over
$L ((t^{1/n}))$
for some finite extension field $L$ of $K$ and some positive integer $n$.

\subsection{Beyond Christol's theorem: generalized power series}

As suggested by Abhyankar \cite{abh}, the situation 
described in the previous section
can be remedied by
allowing certain ``generalized power series''; these were in fact first
introduced by Hahn \cite{hahn} in 1907. We will define these more precisely
in Section~\ref{subsec:genps}; for now, think of
a generalized power series as a series $\sum_{i \in I} x_i t^i$
where the index set $I$ is a well-ordered subset of the rationals
(i.e., a subset containing no infinite decreasing sequence).
For example, in the ring of generalized power series over $\FF_p$,
Chevalley's polynomial has the roots
\[
x = c + t^{-1/p} + t^{-1/p^2} + \cdots
\]
for $c=0,1, \dots, p-1$.

Denote the field of generalized power series over $\FF_q$
by $\FF_q((t^\QQ))$.
Then it turns out that $\FF_q((t^\QQ))$ is algebraically closed, and
one can explicitly characterize those of its elements which are the roots
of polynomials over $\FF_q((t))$ \cite{me}.
One then may ask whether one can, in the vein of Christol, 
give an automata-theoretic characterization
of the elements of $\FF_q((t^\QQ))$ which are roots of monic polynomials over
$\FF_q(t)$. 

In this paper, we give such an automata-theoretic characterization.
(The characterization appeared previously
in the unpublished preprint \cite{me3}; this paper is an updated 
and expanded version of that one.) In the process,
we characterize well-ordered
sets of nonnegative rational numbers with terminating
base $b$ ($b>1$ an integer) which are generated by a finite automaton,
and describe some techniques that may be useful for computing in
the algebraic closure of $\FF_q(t)$, such as an analogue of
Newton's algorithm. (One thing we do not do is give an independent derivation
of Christol's theorem; the new results here are
essentially orthogonal to that result.)
Whether one can use automata in practice to
perform some sort of ``interval arithmetic'' is an intriguing
question about which we will not say anything conclusive,
though we do make a few speculative comments in Section~\ref{subsec:algor}.

\subsection{Structure of the paper}

To conclude this introduction, we describe the contents of 
the remaining chapters of the paper. 

In Chapter~\ref{sec:auto}, we collect some relevant background
material on deterministic finite automata (Section~\ref{subsec:DFA}),
nondeterministic finite automata (Section~\ref{subsec:NFA}),
and the relationship between automata and base $b$ expansions of
rational numbers (Section~\ref{subsec:base}).

In Chapter~\ref{sec:alg}, we collect some relevant background 
material on generalized power series (Section~\ref{subsec:genps}),
algebraic elements of field extensions (Section~\ref{subsec:algebraic}),
and additive polynomials (Section~\ref{subsec:additive}).

In Chapter~\ref{sec:genps-aut}, we state our main theorem
relating generalized power series algebraic over $\FF_q(t)$ with
automata, to be proved later in the paper.
We formulate the theorem and note some
corollaries (Section~\ref{subsec:mainthm}), then refine
the statement by checking its compatibility with ``decimation''
of a power series (Section~\ref{subsec:decimation}).

In Chapter~\ref{sec:abstract}, we give one complete proof of
the main theorem, which in one direction 
relies on a certain amount of sophisticated algebraic
machinery. We give a fairly direct proof that
automatic generalized power series are algebraic
(Section~\ref{subsec:auttoalg}), then
give a proof of the reverse implication by specializing the results of
\cite{me} (Section~\ref{subsec:algtoaut1}); the dependence on
\cite{me} is the source of the reliance on algebraic tools,
such as Artin-Schreier theory.

In Chapter~\ref{sec:valued}, we collect more results about fields
with a valuation, specifically in the case of positive characteristic.
We recall basic properties of ``twisted polynomials''
(Section~\ref{subsec:twisted}) and Newton polygons
(Section~\ref{subsec:np}), give a basic form of Hensel's lemma
on splitting polynomials (Section~\ref{subsec:slope-split}),
and adapt this result to twisted polynomials
(Section~\ref{subsec:slope-twisted}). 

In Chapter~\ref{sec:concrete}, we give a second proof of the
reverse implication of the main theorem, replacing the algebraic
methods of the previous chapter with more explicit considerations
of automata. To do this, we analyze the transition graphs of automata
which give rise to generalized power series
(Section~\ref{subsec:transition}), show that the class of
automatic generalized power series is closed under addition
and multiplication  (Section~\ref{subsec:arithmetic}),
and exhibit a positive-characteristic analogue of Newton's
iteration  (Section~\ref{subsec:Newton}).

In Chapter~\ref{sec:further}, we raise some further questions about
the algorithmics of automatic generalized power series
(Section~\ref{subsec:algor}), and about a potential generalization
of the multivariate analogue of Christol's theorem
(Section~\ref{subsec:multi}).

\subsection*{Acknowledgments}
Thanks to Bjorn Poonen for bringing the work of Christol to the author's
attention, to Eric Rains for some intriguing suggestions
concerning efficient 
representations, to Richard Stanley for pointing out the term
``cactus'', and to George Bergman for helpful 
discussions. The author was supported previously
by an NSF Postdoctoral Fellowship and currently by NSF grant
DMS-0400727.

\section{Automata}
\label{sec:auto}

In this chapter, we recall notions and fix notation and terminology
regarding finite automata. We take as our reference \cite[Chapter~4]{as}.
We note in passing that a sufficiently diligent reader should be able to
reproduce the proofs of all cited results in this chapter.

\subsection{Deterministic automata}
\label{subsec:DFA}

\begin{defn}
A \emph{deterministic finite automaton}, or \emph{DFA} for short, 
is
a tuple $M = (Q, \Sigma, \delta, q_0, F)$, where
\begin{itemize}
\item $Q$ is a finite set (the \emph{states});
\item $\Sigma$ is another set (the \emph{input alphabet});
\item $\delta$ is a function from $Q \times \Sigma$ to $Q$ 
(the \emph{transition function});
\item $q_0 \in Q$ is a state (the \emph{initial state});
\item $F$ is a subset of $Q$ (the \emph{accepting states}).
\end{itemize}
\end{defn}

\begin{defn}
Let $\Sigma^*$ denote the set of finite sequences consisting of
elements of $\Sigma$; we will refer to elements of $\Sigma$ as \emph{characters}
and elements of $\Sigma^*$ as \emph{strings}.
We identify elements of $\Sigma$ with one-element strings,
and denote concatenation of strings by
juxtaposition: that is, if $s$ and $t$ are strings, then $st$ is the
string composed of the elements of $s$ followed by the elements of $t$.
We define a \emph{language} (over $\Sigma$) to be any subset of $\Sigma^*$.
\end{defn}

It is sometimes convenient to represent a DFA using a transition graph.
\begin{defn}
Given a DFA $M = (Q, \Sigma, \delta, q_0, F)$, the 
\emph{transition graph} of $M$ is the edge-labeled directed graph (possibly
with loops) on the vertex set $Q$, with an edge from $q \in Q$ to
$q' \in Q$ labeled by $s \in \Sigma$ if $\delta(q,s) = q'$. The
transition graph also comes equipped with a distinguished vertex
corresponding to $q_0$, and a distinguished subset of the
vertex set corresponding to $F$; from these data, one can recover $M$
from its transition graph.
\end{defn}

One can also imagine a DFA as a machine with a keyboard containing
the elements of $\Sigma$, which can be at any time in any of the states.
When one presses a key, the machine transitions to a new state by applying
$\delta$ to the current state and the key pressed.
One can then extend
the transition function to strings by pressing the corresponding keys in
sequence. Formally, we extend $\delta$ to a function $\delta^*: 
Q \times \Sigma^*
\to Q$ by the rules 
\[
\delta^*(q, \emptyset) = q, \quad
\delta^*(q, xa) = \delta(\delta^*(q,x), a) 
\qquad (q \in Q, x \in \Sigma^*, a \in \Sigma).
\]
\begin{defn}
We say that $M$
\emph{accepts} a string $x \in \Sigma^*$ if $\delta^*(q_0, x) \in 
F$, and otherwise say it \emph{rejects} $x$. The set of strings accepted by
$M$ is called the \emph{language accepted by $M$} and denoted $L(M)$.
A language is said to be \emph{regular} if it is accepted by some DFA.
\end{defn}

\begin{lemma} \label{lem:basic}
\begin{enumerate}
\item[(a)] The collection of regular languages is closed under complement,
finite union, and finite intersection. Also, any language consisting of
a single string is regular.
\item[(b)] The collection of regular languages is closed under reversal
(the operation on strings taking $s_1\cdots s_n$ to $s_n\cdots s_1$).
\item[(c)] A language is regular if and only if it is generated by some
regular expression (see \cite[1.3]{as} for a definition).
\end{enumerate}
\end{lemma}
\begin{proof}
(a) is straightforward, (b) is \cite[Corollary~4.3.5]{as},
and (c) is Kleene's theorem \cite[Theorem~4.1.5]{as}.
\end{proof}

The Myhill-Nerode theorem \cite[Theorem~4.1.8]{as} gives an intrinsic
characterization of regular languages, without reference to an
``auxiliary'' automaton.
\begin{defn}
Given a language $L$ over $\Sigma$, define the equivalence relation
$\sim_L$ on $\Sigma^*$ as follows: 
$x \sim_L y$ if and only if
for all $z \in \Sigma^*$, $xz \in L$ if and only if $yz \in L$.
\end{defn}
\begin{lemma}[Myhill-Nerode theorem] \label{lem:mn}
The language $L$ is regular if and only if $\Sigma^*$ has only finitely
many equivalence classes under $\sim_L$.
\end{lemma}
Moreover, if $L$ is regular, then the DFA in which:
\begin{itemize}
\item $Q$ is the set of equivalence classes under $\sim_L$;
\item $\delta$, applied to the class of some $x \in \Sigma^*$ and some
$s \in \Sigma$, returns the class of $xs$;
\item $q_0$ is the class of the empty string;
\item $F$ is the set of classes of elements of $L$;
\end{itemize}
generates $L$ and has fewer states than any nonisomorphic DFA which also
generates $L$ \cite[Corollary~4.1.9]{as}.

It will also be convenient to permit automata to make non-binary
decisions about strings.
\begin{defn}
A \emph{deterministic finite automaton with output}, 
or \emph{DFAO} for short, 
is
a tuple $M = (Q, \Sigma, \delta, q_0, \Delta, \tau)$, where
\begin{itemize}
\item $Q$ is a finite set (the \emph{states});
\item $\Sigma$ is another set (the \emph{input alphabet});
\item $\delta$ is a function from $Q \times \Sigma$ to $Q$ 
(the \emph{transition function});
\item $q_0 \in Q$ is a state (the \emph{initial state});
\item $\Delta$ is a finite set (the \emph{output alphabet});
\item $\tau$ is a function from $Q$ 
to $\Delta$ (the {output function}).
\end{itemize}
A DFAO $M$ gives rise to a function $f_M: \Sigma^* \to \Delta$
by setting $f_M(w) = \tau(\delta^*(q_0,w))$. Any function $f:
\Sigma^* \to \Delta$ equal to $f_M$ for some DFAO $M$ is called a
\emph{finite-state function}; note that $f$ is a finite-state function
if and only if $f^{-1}(d)$ is a regular language for each $d \in \Delta$
\cite[Theorems~4.3.1 and~4.3.2]{as}.
\end{defn}
We will identify each DFA with the DFAO with output alphabet $\{0,1\}$
which outputs 1 on a string if the original DFA accepts the string and
0 otherwise.

It will also be useful to have devices that can operate on the class
of regular languages.
\begin{defn}
A \emph{finite-state transducer} is a tuple
$T = (Q, \Sigma, \delta, q_0, \Delta, \lambda)$, where
\begin{itemize}
\item $Q$ is a finite set (the \emph{states});
\item $\Sigma$ is another set (the \emph{input alphabet});
\item $\delta$ is a function from $Q \times \Sigma$ to $Q$ 
(the \emph{transition function});
\item $q_0 \in Q$ is a state (the \emph{initial state});
\item $\Delta$ is a finite set (the \emph{output alphabet});
\item $\lambda$ is a function from $Q \times \Sigma$ 
to $\Delta^*$ (the \emph{output function}).
\end{itemize}
If the output of $\lambda$ is always a string of length $k$, we say the
transducer $T$ is \emph{$k$-uniform}. 
\end{defn}

A transducer $T$ gives rise
to a function $f_T: \Sigma^* \to \Delta^*$ as follows: given a
string $w = s_1 \cdots s_r \in \Sigma^*$ (with each $s_i \in \Sigma$),
put $q_i = \delta^*(s_1 \cdots s_i)$ for $i=1, \dots, r$, and 
define
\[
f_T(w) = \lambda(q_0, s_1) \lambda(q_1, s_2) \cdots \lambda(q_{r-1}, s_r).
\]
That is, feed $w$ into the transducer and at each step, use the
current state and the next transition to produce a piece of output,
then string together the outputs.
For $L \subseteq \Sigma^*$ and $L' \subseteq \Delta^*$ languages, 
we write
\begin{align*}
f_T(L) &= \{f_T(w): w \in L \} \\
f_T^{-1}(L') &= \{w \in \Sigma^*: f_T(w) \in L'\}.
\end{align*}
Then one has the following result \cite[Theorems~4.3.6 and~4.3.8]{as}.
\begin{lemma} \label{L:transducer}
Let $T$ be a finite-state transducer. If $L \subseteq \Sigma^*$ is a regular
language, then $f_T(L)$ is also regular; if $L' \subseteq \Delta^*$
is a regular language, then $f_T^{-1}(L')$ is also regular.
\end{lemma}

\subsection{Nondeterministic automata and multiplicities}
\label{subsec:NFA}

Although they do not expand the boundaries of the theory, 
it will be useful in practice to allow so-called ``nondeterministic
automata''.
\begin{defn}
A \emph{nondeterministic finite automaton}, or \emph{NFA} for short,
is a tuple $M = (Q, \Sigma, \delta, q_0, F)$, where
\begin{itemize}
\item $Q$ is a finite set (the \emph{states});
\item $\Sigma$ is another set (the \emph{input alphabet});
\item $\delta$ is a function from $Q \times \Sigma$ to the power set
of $Q$ 
(the \emph{transition function});
\item $q_0 \in Q$ is a state (the \emph{initial state});
\item $F$ is a subset of $Q$ (the \emph{accepting states}).
\end{itemize}
For $M$ an NFA and $w = s_1 \dots s_n \in \Sigma^*$ 
(with $s_i \in \Sigma$ for $i=1,\dots,n$), define an
\emph{accepting path} for $w$ to be a sequence of states
$q_1, \dots, q_n \in Q$ such that
$q_i \in \delta(q_{i-1},s_i)$ for $i=1,\dots, n$ and $q_n \in F$.
Define the \emph{language accepted by $M$} as the set of 
strings $w \in \Sigma^*$ 
for which there exists an accepting path.
\end{defn}
Informally, an NFA is a machine which may make a choice of how to transition
based on the current state and key pressed, or may not be able to make any
transition at all. It accepts a string if there is some way it can transition
from the initial state into an accepting state 
via the corresponding key presses.

Every DFA can be viewed as an NFA, and the language accepted is the same
under both interpretations; hence every language accepted by some DFA
is also accepted by some NFA. The converse is also true 
\cite[Theorem~4.1.3]{as}, so for theoretical purposes, it is typically 
sufficient to work with the conceptually simpler DFAs. However, the conversion
from an NFA of $n$ states may produce a DFA with as many as $2^n$ states,
so in practice this is not usually a good idea.

We will need the following quantitative variant of \cite[Theorem~4.1.3]{as}.
\begin{lemma} \label{lem:quant}
Fix a positive integer $n$.
Let $M = (Q, \Sigma, \delta, q_0, F)$ be an NFA, and let
$f: \Sigma^* \to \ZZ/n\ZZ$ be the function that assigns
to $w \in \Sigma^*$ the number of accepting paths for $w$ in $M$, reduced
modulo $n$. Then $f$ is a finite-state function.
\end{lemma}
\begin{proof}
We construct a DFAO
$M' = (Q', \Sigma', \delta', q'_0, \Delta', \tau')$ 
with the property that $f = f_{M'}$, as follows.
Let $Q'$ denote the set of functions from $Q$ to $\ZZ/n\ZZ$,
and put $\Sigma' = \Sigma$.
Define the function $\delta': Q' \times \Sigma \to Q'$ as follows:
given a function $g: Q \to \ZZ/n\ZZ$ and an element $s \in \Sigma$, let
$\delta'(g,s): Q \to \ZZ/n\ZZ$ be the function given by
\[
\delta'(g, s)(q) = \sum_{q_1 \in Q: \delta(q_1,s) = q} g(q_1). 
\]
Let $q'_0: Q \to \ZZ/n \ZZ$ be the function
carrying $q_0$ to 1 and all other states
to 0.
Put $\Delta' = \ZZ/n\ZZ$, and let $\tau': Q' \to \ZZ/n\ZZ$ be the 
function given by
\[
\tau'(g) = \sum_{q \in F} g(q).
\]
Then $M'$ has the desired properties.
\end{proof}

\subsection{Base expansions and automatic functions}
\label{subsec:base}

In this section, we make precise the notion of ``a function on $\QQ$ computable
by a finite automaton'', and ultimately relate it
to the notion of an ``automatic sequence'' from \cite[Chapter~5]{as}.
To do this, we need to fix a way to input rational numbers into an automaton,
by choosing some conventions about base expansions.

Let $b>1$ be a fixed positive integer. All automata in this section will have
input alphabet $\Sigma = \Sigma_b = \{0, 1, \dots, b-1, .\}$, 
which we identify with the base $b$
digits and radix point. 
\begin{defn}
A string $s = s_1 \dots s_n \in \Sigma^*$ is 
said to be a \emph{valid base $b$ expansion}
if $s_1 \neq 0$, $s_n \neq 0$, and exactly one of $s_1,\dots,s_n$
is equal to the radix point. If $s$ is a valid base $b$ expansion
and $s_k$ is the radix point,
 then we define the
\emph{value} of $s$ to be
\[
v(s) = \sum_{i=1}^{k-1} s_i b^{k-1-i} + \sum_{i=k+1}^n s_i b^{k-i}.
\]
It is clear that no two valid strings have the same value; we may thus
unambiguously define $s$ to be the \emph{base $b$ expansion of $v(s)$}.
Let $S_b$ be the set of nonnegative $b$-adic rationals, i.e., numbers of the
form $m/b^n$ for some nonnegative integers $m,n$; it is also clear
that the set of values of valid base $b$ expansions is precisely $S_b$.
For $v \in S_b$, write $s(v)$ for the base $b$ expansion of $v$.
\end{defn}

\begin{lemma} \label{lem:valid}
The set of valid base $b$ expansions is a regular language.
\end{lemma}
\begin{proof}
The language $L_1$ of strings with no leading zero is regular by
Lemma~\ref{lem:mn}: the equivalence classes under $\sim_{L_1}$ consist
of the empty string, all nonempty strings in $L_1$, and all nonempty 
strings not in $L_1$. The language $L_2$ of strings with no trailing
zero is also regular: the equivalence classes under $\sim_{L_2}$
consist of all strings in $L_2$, and all
strings not in $L_2$. (One could also apply Lemma~\ref{lem:basic}(b)
to $L_1$ to show that $L_2$ is regular, or vice versa.) The language $L_3$
of strings with exactly one radix point is also regular: the equivalence
classes under $\sim_{L_3}$ consist of all strings with zero points,
all strings with one point, and all strings with more than one point. Hence
$L_1 \cap L_2 \cap L_3$ is regular by Lemma~\ref{lem:basic}(a), as desired.
\end{proof}

The ``real world'' convention for base $b$ expansions is a bit more
complicated than what we are using: normally, one omits the radix point
when there are no digits after it, one adds a leading zero in front
of the radix point when there are no digits before it, and one
represents 0 with a single zero rather than a bare radix point (or the
empty string). This will not
change anything essential, thanks to the following lemma.
\begin{lemma} \label{lem:realworld}
Let $S$ be a set of nonnegative $b$-adic rationals. Then the set of 
expansions of $S$, under our convention, is a regular language if and only
if the set of ``real world'' expansions of $S$ is a regular language.
\end{lemma}
\begin{proof}
We may as well assume for simplicity that $0 \notin S$ since any
singleton language is regular. Put
\[
S_1 = S \cap (0,1), \qquad S_2 = S \cap \ZZ, \qquad S_3 = S \setminus 
(S_1 \cup S_2).
\]
Then the expansions of $S$, under either convention, form a regular language
if and only if the same is true of $S_1, S_2, S_3$. Namely, 
under our convention, the language 
of strings with no digits before the radix point
and the language of strings with no digits after the radix point are regular.
Under the ``real world'' convention, the language
 of strings with no radix point
and the language of strings with a single 0 before the radix point 
are regular.

The expansions of $S_3$ are the same in both cases, so we can ignore them.
For $S_1$, note that the language of its real world expansions is regular if 
and only if the language of the reverses of those strings 
is regular (Lemma~\ref{lem:basic}(b)),
if and only if the languages of those reverses with
a radix point added in front are regular (clear),
if and only if the languages of the reverses of 
those (which are the expansions under our
convention) are regular.
For $S_2$, note that the language of real world expansions is regular if 
and only if the language of those strings
with the initial zeroes removed is regular.
\end{proof}

\begin{defn}
Let $M$ be a DFAO with input alphabet $\Sigma_b$. We say a state $q \in Q$
is \emph{preradix} (resp.\ \emph{postradix}) if there exists a valid 
base $b$ expansion $s = s_1\cdots s_n$ with $s_k$ equal to the radix point
such that, if we set $q_i = \delta(q_{i-1}, s_i)$, then $q = q_i$
for some $i < k$ (resp.\ for some $i \geq k$). That is, when tracing through
the transitions produced by $s$, $q$ appears before (resp.\ after) the
transition producing the radix point. Note that if the language
accepted by $M$ consists only of valid base $b$ expansions, then
no state can be both preradix and postradix, or else $M$ would accept
some string containing more than one radix point.
\end{defn}

\begin{defn}
Let $\Delta$ be a finite set. A function $f: S_b \to \Delta$ is
\emph{$b$-automatic} if there is a DFAO $M$ with input alphabet $\Sigma$
and output alphabet $\Delta$ such that for any $v \in S_b$, 
$f(v) = f_M(s(v))$. By Lemma~\ref{lem:valid}, it is equivalent to require
that for some symbol $\star \notin \Delta$, there is a DFAO $M$ with
input alphabet $\Sigma$ and output alphabet $\Delta \cup \{\star\}$ such that
\[
f_M(s) = \begin{cases} f(v(s)) & \mbox{$s$ is a valid base $b$ expansion} \\
\star & \mbox{otherwise.} 
\end{cases}
\]
We say a subset $S$ of $S_b$ is \emph{$b$-regular} if its characteristic
function
\[
\chi_S(s) = \begin{cases} 1 & s \in S \\ 0 & s \notin S \end{cases}
\]
is $b$-automatic; then a function $f: S_b \to \Delta$ is $b$-automatic
if and only if $f^{-1}(d)$ is $b$-regular for each $d \in \Delta$.
\end{defn}

\begin{lemma} \label{L:dilation}
Let $S \subseteq S_b$ be a subset. Then for any $r \in \NN$
and any $s \in S_b$, $S$ is $b$-regular if and only if
\[
rS + s = \{rx + s: x \in S\}
\]
is $b$-regular.
\end{lemma}
\begin{proof}
As in \cite[Lemmas~4.3.9 and~4.3.11]{as}, 
one can construct a finite-state transducer that performs the
operation $x \mapsto rx + s$ on valid base $b$ expansions read
from right to left, by simply transcribing the usual hand
calculation. (Remember that reversing the strings of a language
preserves regularity by Lemma~\ref{lem:basic}, so there is no
harm in reading base $b$ expansions backwards.)
Lemma~\ref{L:transducer} then yields the desired result.
\end{proof}

We conclude this section by noting the relationship with the notion of
``automatic sequences'' from \cite[Chapter~5]{as}. In 
\cite[Definition~5.1.1]{as}, a sequence
$\{a_l\}_{l=0}^\infty$ over $\Delta$
is said to be \emph{$b$-automatic} if there is a DFAO $M$ with input alphabet
$\{0, \dots, b-1\}$ and output alphabet $\Delta$ such that for any string $s
= s_1\cdots s_n$, if we put $v(s) = \sum_{i-1}^n s_1 b^{n-1-i}$, then
$a_{v(s)} = f_M(s)$. Note that this means $M$ must evaluate correctly even
on strings with leading zeroes, but by \cite[Theorem~5.2.1]{as}, it is
equivalent to require that there exists such an $M$ only having the
property that $a_{v(s)} = f_M(s)$ when $s_1 \neq 0$.
It follows readily that $\{a_l\}_{l=0}^\infty$ is $b$-automatic if and only
if for some symbol $\star \notin \Delta$, the function $f: S_b \to \Delta
\cup \{ \star\}$ defined by
\[
f(x) = \begin{cases} a_{x} & x \in \ZZ \\
\star & \mbox{otherwise}
\end{cases}
\]
is $b$-automatic.

\section{Algebraic preliminaries}
\label{sec:alg}

In this chapter, we recall the algebraic machinery that will go into
the formulation of Theorem~\ref{thm:main}.

\subsection{Generalized power series}
\label{subsec:genps}

Let $R$ be an arbitrary ring.
Then the ring of ordinary power series over $R$ can be identified with
the ring of functions from $\ZZ_{\geq 0}$ to $R$, with addition given termwise
and multiplication given by convolution
\[
(fg)(k) = \sum_{i+j=k} f(i) g(j);
\]
the latter makes sense because for any fixed $k \in \ZZ_{\geq 0}$,
there are only finitely many pairs $(i,j) \in \ZZ_{\geq 0}^2$ such that
$i+j=k$. In order to generalize this construction to index sets other
than $\ZZ_{\geq 0}$, we will have to restrict the nonzero values of
the functions so that computing $fg$ involves adding only finitely many
nonzero elements of $R$. The recipe for doing this dates back to
Hahn \cite{hahn} (although the term ``Mal'cev-Neumann series'' for
an object of the type we describe is prevalent), and we recall it now; see also
\cite[Chapter~13]{passman}.

\begin{defn}
Let $G$ be a totally ordered abelian group (written additively) with
identity element $0$; that is, $G$ is an abelian group equipped with
a binary relation $>$ such that for all $a,b,c \in G$,
\begin{gather*}
a \not> a \\
a \not> b, b \not> a \Rightarrow a = b \\
a > b, b > c \Rightarrow a > c \\
a > b \Leftrightarrow a+c > b+c.
\end{gather*}
Let $P$ be the set of $a \in G$ for which $a>0$; $P$ is called the
\emph{positive cone} of $G$.
\end{defn}
\begin{lemma}\label{lem:wo}
Let $S$ be a subset of $G$. Then the
following two conditions are equivalent.
\begin{enumerate}
\item[(a)] Every nonempty subset of $S$ has a minimal element.
\item[(b)] There is no infinite decreasing sequence $s_1 > s_2 > \cdots$
within $S$.
\end{enumerate}
\end{lemma}
\begin{proof}
If (a) holds but (b) did not, then the set $\{s_1, s_2, \dots\}$ would not
have a minimal element, a contradiction. Hence (a) implies (b). Conversely,
if $T$ were a subset of $S$ with no smallest element, then for any
$s_i \in T$, we could choose $s_{i+1} \in T$ with $s_i > s_{i+1}$, thus
forming an infinite decreasing sequence. Hence (b) implies (a).
\end{proof}
\begin{defn}
A subset $S$ of $G$ is \emph{well-ordered} if it satisfies either of the
equivalent conditions of Lemma~\ref{lem:wo}. (Those who prefer to avoid
assuming the axiom of choice should take (a) to be the definition, as
the implication (b)$\implies$(a) requires choice.)
\end{defn}

For $S_1, \dots, S_n \subseteq G$, write $S_1 + \cdots + S_n$ for the
set of elements of $G$ of the form $s_1 + \cdots + s_n$ for $s_i \in S_i$;
in case $S_1 = \cdots = S_n$, we abbreviate this notation to
$S^{+n}$. (In \cite{passman} the notation
$nS$ is used instead, but we have already defined this as the dilation
of $S$ by the factor $n$.)
Then one can easily verify the following (or see
\cite[Lemmas~13.2.9 and~13.2.10]{passman}).
\begin{lemma} \label{L:passman}
\begin{enumerate}
\item[(i)] 
If $S_1, \dots, S_n$ are well-ordered subsets of $G$, then 
$S_1 + \cdots + S_n$ is well-ordered.
\item[(ii)]
If $S_1, \cdots, S_n$ are well-ordered subsets of $G$, then
for any $x \in G$, the number of $n$-tuples $(s_1, \dots, s_n) \in
S_1 \times \cdots \times S_n$ such that $s_1 + \cdots + s_n = x$ is finite.
\item[(iii)]
If $S$ is a well-ordered subset of $P$, then
$\tilde{S} = \cup_{n=1}^\infty S^{+n}$
also is well-ordered; moreover, $\cap_{n=1}^\infty \tilde{S}^{+n} = \emptyset$.
\end{enumerate}
\end{lemma}

\begin{defn}
Given a function $f: G \to R$, the 
\emph{support} of $f$ is the set of $g \in G$ such that $f(g) \neq 0$.
A \emph{generalized Laurent series} over $R$ with exponents in $G$
is a function $f: G \to R$ whose support is well-ordered;
if the support is contained in $P \cup \{0\}$, we call $f$ a \emph{generalized
power series}. We typically represent the generalized Laurent series $f$
in series notation $\sum_i f(i) t^i$, and write $R\llbracket t^G \rrbracket$
and $R((t^G))$ for the sets of generalized power series and generalized
Laurent series, respectively, over $R$ with exponents in $G$.
\end{defn}

Thanks to Lemma~\ref{L:passman},
the termwise sum and convolution product are well-defined binary
operations on $R \llbracket t^G \rrbracket$ and $R ((t^G))$, which form 
rings under the operations. 
A nonzero element of
$R((t^G))$ if and only if its first nonzero coefficient is a unit
\cite[Theorem~13.2.11]{passman};
in particular, if $R$ is a field, then
so is $R((t^G))$.

\subsection{Algebraic elements of fields}
\label{subsec:algebraic}

In this section, we recall the definition of algebraicity of an element
of one field over a subfield, and then review some criteria for
algebraicity. Nothing in this section is even remotely original, as can be confirmed by any sufficiently detailed abstract algebra textbook.

\begin{defn}
Let $K \subseteq L$ be fields. Then $\alpha \in L$ is said to be
\emph{algebraic over $K$} if there exists a nonzero 
polynomial $P(x) \in K[x]$ over $K$ such that $P(\alpha) = 0$.
We say $L$ is \emph{algebraic over $K$} if every element of $L$ is
algebraic over $K$. 
\end{defn}
\begin{lemma}
Let $K \subseteq L$ be fields. Then $\alpha \in L$ is algebraic if and only
if $\alpha$ is contained in a subring of $L$ containing $K$
which is finite dimensional
as a $K$-vector space.
\end{lemma}
\begin{proof}
We may as well assume $\alpha \neq 0$, as otherwise both assertions are clear.
If $\alpha$ is contained in a subring $R$ of $L$ which has
finite dimension $m$ as a $K$-vector space, then $1, \alpha, \dots, \alpha^m$
must be linearly dependent over $K$, yielding a polynomial over $K$
with $\alpha$ as a root. Conversely, if $P(\alpha) = 0$ for some polynomial
$P(x) \in K[x]$, we may take $P(x) = c_0 + c_1 x + \cdots + c_n x^n$
with $c_0, c_n \neq 0$. In that case, 
\[
\{a_0 + a_1 \alpha + \cdots + a_{n-1} \alpha^{n-1}: a_0,\dots,a_{n-1} \in K\}
\]
is a subring of $L$ containing $K$ and $\alpha$, of dimension at most $n$
as a vector space over $K$.
\end{proof}
\begin{cor}
Let $K \subseteq L$ be fields. If $\dim_K L < \infty$, then $L$ is
algebraic over $K$.
\end{cor}

\begin{lemma} \label{L:alg transitive}
Let $K \subseteq L \subseteq M$ be fields, with $L$ algebraic
over $K$. For any $\alpha \in M$, $\alpha$ is algebraic over $K$
if and only if it is algebraic over $L$.
\end{lemma}
\begin{proof}
Clearly if $\alpha$ is the root of a polynomial with coefficients
in $K$, that same polynomial has coefficients in $L$. Conversely,
suppose $\alpha$ is algebraic over $L$; it is then contained
in a subring $R$ of $M$ which is finite dimensional over $L$. That subring
is generated over $L$ by finitely many elements, each of which is
algebraic over $K$ and hence lies in a subring $R_i$ of $M$ which is
finite dimensional over $K$. Taking the ring generated by the $R_i$
gives a subring of $M$ which is finite dimensional over $K$ and
which contains $\alpha$. Hence $\alpha$ is algebraic over $K$.
\end{proof}

\begin{lemma} \label{lem:subring}
Let $K \subseteq L$ be fields. If $\alpha,\beta \in L$ are algebraic over $K$,
then so are $\alpha+\beta$ and $\alpha \beta$. If $\alpha \neq 0$, then
moreover $1/\alpha$ is algebraic over $K$.
\end{lemma}
\begin{proof}
Suppose that $\alpha, \beta \in L$ are algebraic over $K$;
we may assume $\alpha, \beta \neq 0$, else everything is clear.
Choose polynomials $P(x) = c_0 + c_1 x + \cdots + c_m x^m$ and
$Q(x) = d_0 + d_1 x + \cdots + d_n x^n$ with $c_0, c_m, d_0, d_n \neq 0$
such that $P(\alpha) = Q(\beta) = 0$. Then 
\[
R = \left\{ \sum_{i=0}^{m-1} \sum_{j=0}^{n-1} a_{ij} \alpha^i \beta^j:
a_{ij} \in K
\right\}
\]
is a subring of $L$ containing $K$, of dimension at most $mn$ as a
vector space over $K$, containing $\alpha + \beta$ and $\alpha \beta$.
Hence both of those are algebraic over $K$. Moreover, 
$1/\alpha = -(c_1 + c_2 \alpha + \cdots + c_m x^{m-1})/c_0$ is contained
in $R$, so it too is algebraic over $K$.
\end{proof}

\begin{defn}
A field $K$ is \emph{algebraically closed} if every polynomial over $K$
has a root, or equivalently, if every polynomial over $K$ splits
completely (factors into linear polynomials). 
It can be shown (using Zorn's lemma)
that every field $K$ is contained in an algebraically closed
field; the elements of such a field 
which are algebraic over $K$ form a field $L$
which is both algebraically closed and algebraic over $K$. Such a field
is called an \emph{algebraic closure} of $K$; it can be shown to be unique
up to noncanonical isomorphism, but we won't need this.
\end{defn}

In practice, we will always consider fields contained in $\FF_q((t^\QQ))$,
and constructing algebraically closed fields containing them is 
straightforward. That is because
if $K$ is an algebraically closed field and $G$ is a divisible group
(i.e., multiplication by any positive integer is a bijection on $G$),
then the field $K((t^G))$ is algebraically closed. 
(The case $G = \QQ$, which is
the only case we need, is treated explicitly in \cite[Proposition~1]{me2};
for the general case and much more, see \cite[Theorem~5]{kap}.)
Moreover, it is easy (and does not require the axiom of choice)
to construct an algebraic closure $\overline{\FF_q}$
of $\FF_q$: order the elements of
$\FF_q$ with 0 coming first, 
then list the monic polynomials over $\FF_q$ in lexicographic
order and successively adjoin roots of them.
Then the field $\overline{\FF_q}((t^\QQ))$ is algebraically closed
and contains $\FF_q((t^\QQ))$.

\subsection{Additive polynomials}
\label{subsec:additive}

In positive characteristic, it is convenient to restrict attention to
a special class of polynomials, the ``additive'' polynomials.
First, we recall a standard recipe (analogous to the construction of
Vandermonde determinants) for producing such polynomials.

\begin{lemma} \label{L:moore}
Let $K$ be a field of characteristic $p>0$.
Given $r_1, \dots, r_n \in K$, the
Moore determinant
\[
\begin{pmatrix}
r_1 & r_2 & \cdots & r_n \\
r_1^p & r_2^p &\cdots & r_n^{p} \\
\vdots & \vdots & \ddots & \vdots \\
r_1^{p^{n-1}} & r_2^{p^{n-1}} & \cdots & r_n^{p^{n-1}}
\end{pmatrix}
\]
vanishes if and only if $r_1, \dots, r_n$ are linearly dependent over
$\FF_p$.
\end{lemma}
\begin{proof}
Viewed as a polynomial in $r_1, \dots, r_n$ over $\FF_p$, the Moore determinant
is divisible by each of the linear forms
$c_1 r_1 + \cdots + c_n r_n$ for $c_1, \dots, c_n \in \FF_p$
not all zero. Up to scalar multiples, there are $p^{n-1} + \cdots + p + 1$
such forms, so the determinant is divisible by the product of these forms.
However, the determinant visibly is a homogeneous polynomial in the
$r_i$ of degree $p^{n-1} + \cdots + p + 1$, so it must be equal to the product
of the linear factors times a constant. The desired result follows.
\end{proof}

\begin{defn}
A polynomial $P(z)$ over a field $K$ of characteristic $p>0$ is
said to be \emph{additive} (or \emph{linearized}) if it has the form
\[
P(z) = c_0 z + c_1 z^p + \cdots + c_n z^{p^n}
\]
for some $c_0, \dots, c_n \in K$.
\end{defn}

\begin{lemma} \label{L:additive}
Let $P(z)$ be a nonzero polynomial over a field $K$ of characteristic $p>0$,
and let $L$ be an algebraic closure of $K$. Then the following conditions
are equivalent.
\begin{enumerate}
\item[(a)]
The polynomial $P(z)$ is additive.
\item[(b)]
The equation $P(y+z) = P(y) + P(z)$ holds as a formal identity of polynomials.
\item[(c)]
The equation $P(y+z) = P(y) + P(z)$ holds for all $y,z \in L$.
\item[(d)]
The roots of $P$ in $L$ form an $\FF_p$-vector space under addition,
all roots occur to the same multiplicity, and that multiplicity is a
power of $p$.
\end{enumerate}
\end{lemma}
\begin{proof}
The implications (a)$\implies$(b)$\implies$(c) are clear, and
(c)$\implies$(b) holds because the field $L$ must be infinite.
We next check that (d)$\implies$(a). Let $V \subset L$ be the set of roots,
and let $p^e$ be the common multiplicity.
Let $Q(z)$ be the Moore determinant of $z^{p^e}, r_1^{p^e}, \cdots, r_m^{p^e}$.
By Lemma~\ref{L:moore}, the roots of $Q$ are precisely the elements of $V$,
and each occurs with multiplicity at least $p^e$.
However, $\deg(Q) = p^{e+m} = \deg(P)$, 
so the multiplicities must be exactly $p^e$, and $P$ must equal $Q$ times
a scalar. Since $Q$ is visibly additive, so is $P$.

It remains to check that (c)$\implies$(d).
Given (c), note that the roots of $P$ in $L$ form an $\FF_p$-vector space
under addition;
also, if $r \in L$ is a root of $P$, then $P(z+r) = P(z)$, so all roots
of $P$ have the same multiplicity.
Let $V$ be the roots of $P$, choose generators $r_1, \dots, r_m$
of $V$ as an $\FF_p$-vector space, 
and let $Q(z)$ be the Moore determinant of $z, r_1, \dots, r_m$.
Then $P(z) = c Q(z)^n$ for some
constant $c$, where $n$ is the common multiplicity of the roots of $P$
(because $Q$ has no repeated roots, by the analysis of the
previous paragraph).
Suppose that $n$ is not a prime power; then the polynomials $(y+z)^n$
and $y^n + z^n$ are not identically equal, because the binomial coefficient
$\binom{n}{p^i}$, for $i$ the largest integer such that $p^i$ divides $n$,
is not divisible by $p$. Thus there exist values
of $y,z$ in $L$ for which $(y+z)^n \neq y^n + z^n$.
Since $L$ is algebraically closed, $Q$ is surjective as a map from $L$
to itself; we can thus choose $y,z \in L$ such that $(Q(y) + Q(z))^n
\neq Q(y)^n + Q(z)^n$. Since $Q$ is additive, this means that
$P(y+z) \neq P(y) + P(z)$, contrary to hypothesis. We conclude that
$n$ must be a prime power. Hence (c)$\implies$(d), and the proof is complete.
\end{proof}

The following observation is sometimes known as ``Ore's lemma'' (as in
\cite[Lemma~12.2.3]{as}).
\begin{lemma} \label{L:ore}
 For $K \subseteq L$ fields of characteristic $p>0$ and $\alpha \in L$,
$\alpha$ is algebraic over $K$ if and only if it is a root of some
additive polynomial over $K$.
\end{lemma}
\begin{proof}
Clearly if $\alpha$ is a root of an additive polynomial over $K$,
then $\alpha$ is algebraic over $K$. Conversely,
if $\alpha$ is algebraic, then $\alpha, \alpha^p, \dots$ cannot all be
linearly independent, so there must be a linear relation of the form
$c_0 \alpha + c_1 \alpha^p + \cdots + c_n \alpha^{p^n} = 0$ with
$c_0, \dots, c_n \in K$ not all zero. 
\end{proof}

Our next lemma generalizes Lemma~\ref{L:ore} to
``semi-linear'' systems of equations.
\begin{lemma} \label{lem:vector}
Let $K \subseteq L$ be fields of characteristic $p>0$, let $A,B$ be
 $n \times n$ matrices with entries in $K$,
at least one of which is invertible, and let $\bw \in K^n$ be any
(column) vector.
Suppose $\bv \in L^n$ is a vector such that
$A\bv^\sigma + B\bv = \bw$, where $\sigma$ denotes the $p$-th power Frobenius
map. Then the entries of $\bv$ are algebraic over $K$.
\end{lemma}
\begin{proof}
Suppose $A$ is invertible. Then for $i=1,2,\dots$, we can write
$\bv^{\sigma^i} = U_i \bv + \bw_i$ for some $n \times n$
matrix $U_i$ over $K$ and some $\bw_i \in K^n$. 
Such vectors span a vector space over
$K$ of dimension at most $n^2 + n$; for some $m$, we can thus find
$c_0, \dots, c_m$ not all zero such that
\begin{equation} \label{eq:lincom}
c_0 \bv + c_1 \bv^\sigma + \cdots + c_m \bv^{\sigma^m} = 0. 
\end{equation}
Apply Lemma~\ref{L:ore} to each component in \eqref{eq:lincom} to deduce
that the entries of $\bv$ are algebraic over $K$.

Now suppose $B$ is invertible. 
There is no harm in enlarging $L$, so we may as well assume that $L$ is
closed under taking $p$-th roots, i.e., $L$ is perfect. Then the map
$\sigma: L \to L$ is a bijection. Let $K'$ be the set of $x \in L$ for
which there exists a nonnegative integer $i$ such that $x^{\sigma^i} \in K$;
then $\sigma: K' \to K'$ is also a bijection, and
each element of $K'$ is algebraic over $K$.

For $i=1,2,\dots$, we can now write
$\bv^{\sigma^{-i}} = U_i \bv + \bw_i$ for some $n \times n$ matrix
$U_i$ over $K'$ and some $\bw_i \in (K')^n$. As above, we conclude that
the entries of $\bv$ are algebraic over $K'$. However, any element
$\alpha \in L$ algebraic over $K'$ is algebraic over $K$: if
$d_0 + d_1 \alpha + \cdots + d_m \alpha^m = 0$ for $d_0, \dots, d_m \in K'$
not all zero,
then we can choose a nonnegative integer $i$ such that
$d_0^{\sigma^i}, \dots, d_m^{\sigma^i}$ belong to $K$, and
$d_0^{\sigma^i} + d_1^{\sigma^i} \alpha^{p^i} + \cdots + 
d_m^{\sigma^i} \alpha^{mp^i} = 0$.
We conclude that the entries of $\bv$ are algebraic over $K$, as desired.
\end{proof}

\section{Generalized power series and automata}
\label{sec:genps-aut}

In this chapter, we state the main theorem (Theorem~\ref{thm:main})
and some related results; its proof (or rather proofs) will occupy
much of the rest of the paper.

\subsection{The main theorem: statement and preliminaries}
\label{subsec:mainthm}

We are now ready to state our generalization of Christol's theorem,
the main theoretical result of this paper.
For context, we first
 state a form of Christol's theorem (compare
\cite{christol}, \cite{ckmr}, and also \cite[Theorem~12.2.5]{as}).
Reminder: $\FF_q(t)$ denotes the field of rational functions over $\FF_q$,
i.e., the field of fractions of the ring of polynomials $\FF_q[t]$.
\begin{theorem} [Christol] \label{T:christol}
Let $q$ be a power of the prime $p$, and
let $\{a_i\}_{i=0}^\infty$ be a sequence over $\FF_q$. Then
the series $\sum_{i=0}^\infty a_i t^i \in \FF_q \llbracket t \rrbracket$ is 
algebraic over $\FF_q(t)$ if and only if the sequence $\{a_i\}_{i=0}^\infty$
is $p$-automatic.
\end{theorem}

We now formulate our generalization of Christol's theorem.
Recall that $S_p$ is the set of numbers 
of the form $m/p^n$, for $m,n$ nonnegative integers.
\begin{defn} \label{D:automatic series}
Let $q$ be a power of the prime $p$, and
let $f: \QQ \to \FF_q$ be a function whose support $S$
is well-ordered.
We say the generalized Laurent series $\sum_i f(i) t^i$ is
\emph{$p$-quasi-automatic} if the following conditions hold.
\begin{enumerate}
\item[(a)]
For some integers $a$ and $b$ with $a > 0$, the set
$aS+b = \{ai+b: i \in S\}$ is contained in $S_p$, i.e., consists of
nonnegative $p$-adic rationals.
\item[(b)]
For some $a,b$ for which (a) holds, the function
$f_{a,b}: S_p \to \FF_q$ given by 
$f_{a,b}(x) = f((x-b)/a)$ is $p$-automatic.
\end{enumerate}
Note that by Lemma~\ref{L:dilation}, if (b) holds for a single choice
of $a,b$ satisfying (a), then (b) holds also for any choice of $a,b$
satisfying (a).
In case
(a) and (b) hold with $a=1,b=0$, we say the series is \emph{$p$-automatic}.
\end{defn}

\begin{theorem} \label{thm:main}
Let $q$ be a power of the prime $p$, and
let $f: \QQ \to \FF_q$ be a function whose support
is well-ordered.
Then the corresponding generalized Laurent series $\sum_i f(i) t^i
\in \FF_q (( t^\QQ ))$ 
is algebraic over $\FF_q(t)$ if and only if it is $p$-quasi-automatic.
\end{theorem}
We will give two proofs of Theorem~\ref{thm:main} in due course.
In both cases, we use Proposition~\ref{P:auttoalg} to deduce the
implication ``automatic implies algebraic''. For the reverse implication
``algebraic implies automatic'', we use Proposition~\ref{P:algtoaut1}
for a conceptual proof and Proposition~\ref{P:algtoaut2}
for a more algorithmic proof. Note, however, that both of the proofs
in this direction rely on Christol's theorem, so we do not obtain
an independent derivation of that result.

\begin{cor}
The generalized Laurent series $\sum_i f(i) t^i \in \FF_q((t^\QQ))$ is
algebraic over $\FF_q(t)$ 
if and only if for each $\alpha \in \FF_q$, the generalized
Laurent series
\[
\sum_{i \in f^{-1}(\alpha)} t^i
\]
is algebraic over $\FF_q(t)$.
\end{cor}

We mention another corollary following \cite[Theorem~12.2.6]{as}. 
\begin{defn}
Given two
generalized Laurent series $x = \sum_i x_i t^i$ and $y = \sum_i y_i t^i$
in $\FF_q((t^\QQ))$, then $\sum_i (x_iy_i) t^i$ is also a generalized
Laurent series; it is called the \emph{Hadamard product} and denoted
$x \odot y$. Then one has the following assertion, which in the
case of ordinary power series is due to Furstenberg \cite{furstenberg}.
\end{defn}
\begin{cor}
If $x,y \in \FF_q((t^\QQ))$ are algebraic over $\FF_q(t)$, then so
is $x \odot y$.
\end{cor}
\begin{proof}
Thanks to Theorem~\ref{thm:main},
this follows from the fact that if $f: \Sigma^* 
\to \Delta_1$ and $g: \Sigma^* \to \Delta_2$ are finite-state functions,
then so is $f \times g: \Sigma^* \to \Delta_1 \times \Delta_2$; 
the proof of the latter is straightforward
(or compare \cite[Theorem~5.4.4]{as}).
\end{proof}

\subsection{Decimation and algebraicity}
\label{subsec:decimation}

Before we attack Theorem~\ref{thm:main} proper, it will be helpful
to know that
the precise choice of $a,b$ in Theorem~\ref{thm:main},
which does not matter on the automatic side
(Definition~\ref{D:automatic series}), also does not matter on the
algebraic side.

\begin{defn}
For $\tau \in \Gal(\FF_q/\FF_p)$, regard $\tau$ as an automorphism of
$\FF_q(t)$ and $\FF_q((t^\QQ))$ by allowing it to act on coefficients.
That is,
\[
\left( \sum_i x_i t^i \right)^\tau = \sum_i x_i^\tau t^i.
\]
Let $\sigma \in \Gal(\FF_q/\FF_p)$ denote the $p$-power Frobenius map;
note that the convention we just introduced means that
$x^p = x^\sigma$ if $x \in \FF_q$, but not if $x \in \FF_q(t)$
or $x \in \FF_q((t^\QQ))$.
\end{defn}

\begin{lemma} \label{L:decimate}
Let $a,b$ be integers with $a>0$. Then
$\sum_i x_i t^i \in \FF_q((t^\QQ))$ 
is algebraic over $\FF_q(t)$ if and only if
$\sum_i x_{ai+b} t^i$ is algebraic over $\FF_q(t)$.
\end{lemma}
\begin{proof}
It suffices to prove the result in the case $a=1$ and in the case $b=0$,
as the general case follows by applying these two in succession.
The case $a=1$ is straightforward: if
$x = \sum_i x_i t^i$ is a root of the polynomial $P(z)$ over $\FF_q(t)$, then
$x' = \sum_i x_{i+b} t^i = \sum_i x_i t^{i-b}$ is a root of the 
polynomial $P(z t^b)$, and vice versa.

As for the case $b=0$, we can further break it down into two cases,
one in which $a=p$, the other in which $a$ is coprime to $p$. We treat
the former case first.
If $x = \sum_i x_i t^i$ is a root of the polynomial $P(z) = \sum c_j z^j$ over
$\FF_q(t)$, then $x' = \sum_i x_{pi} t^i = \sum_i x_i t^{i/p}$ is 
a root of the polynomial
\[
\sum c_j^\sigma z^{pj}
\]
over $\FF_q(t)$.
Conversely, if $x'$ is a root of the polynomial $Q(z) = \sum d_j z^j$
over $\FF_q(t)$, then $x$ is a root of the polynomial
\[
\sum (d_j^p)^{\sigma^{-1}} z^j
\]
over $\FF_q(t)$.

Now suppose that $b=0$ and $a$ is coprime to $p$. 
Let $\tau: \FF_q((t^\QQ)) \to \FF_q((t^\QQ))$ denote the automorphism
$\sum x_i t^i \mapsto \sum x_i t^{ai}$; then $\tau$ also acts on
$\FF_q(t)$. 
If $x = \sum_i x_i t^i$ is a root of the polynomial $P(z) = \sum c_j z^j$ over
$\FF_q(t)$, then $x' = \sum_i x_{ai} t^i = \sum_i x_i t^{i/a}$ is 
a root of the polynomial
\[
\sum c_j^{\tau^{-1}} z^{j}
\]
over $\FF_q(t^{1/a})$; since $\FF_q(t^{1/a})$ is finite dimensional over
$\FF_q(t)$, $x'$ is algebraic over $\FF_q(t)$ by 
Lemma~\ref{L:alg transitive}.
Conversely, if $x'$ is a root of the polynomial $Q(z) = \sum d_j z^j$
over $\FF_q(t)$, then $x$ is a root of the polynomial
\[
\sum c_j^\tau z^j
\]
over $\FF_q(t)$.

We have now proved the statement of the lemma in case $a=1$ and $b$
is arbitrary, in case $a=p$ and $b=0$, and in case $a$ is coprime to
$p$ and $b=0$. As noted above, these three cases together imply the
desired result.
\end{proof}

\section{Proof of the main theorem: abstract approach}
\label{sec:abstract}

In this chapter, we give a proof of Theorem~\ref{thm:main}.
While the proof in the ``automatic implies algebraic'' direction
is fairly explicit, the proof in the reverse direction relies on
the results of \cite{me}, and hence is fairly conceptual.
We will give a more explicit approach to the reverse direction
in the next chapter.

\subsection{Automatic implies algebraic}
\label{subsec:auttoalg}

In this section, we establish the ``automatic implies algebraic''
direction of Theorem~\ref{thm:main}. The proof is a slight
modification of the usual argument used to prove the corresponding
direction of Christol's theorem
(as in \cite[Theorem~12.2.5]{as}).
(Note that this direction of Theorem~\ref{thm:main} will be invoked in
both proofs of the reverse direction.)

\begin{lemma} \label{L:auttoalg}
Let $p$ be a prime number, and let $S$ be a $p$-regular subset of $S_p$.
Then $\sum_{i \in S} t^i \in \FF_p \llbracket t^{\QQ} \rrbracket$ is algebraic
over $\FF_p(t)$.
\end{lemma}
\begin{proof}
Let $L$ be the language of strings of the form
$s(v)$ for $v \in S$, and
let $M$ be a DFA which accepts $L$.

For $n$ a nonnegative integer, let $s'(n)$ be the base $p$ expansion of 
$n$ minus the final radix point.
For each preradix state $q \in Q$, let $T_q$ be the set of nonnegative
integers $n$ such that $\delta^*(q_0, s'(n)) = q$,
put $f(q) = \sum_{i \in T_q} t^i$,
and let $U_q$ be the set of pairs $(q',d) \in Q \times \{0,
\dots, p-1\}$ such that $\delta(q',d) = q$. (Note that this forces $q'$
to be preradix.) Then if $q \neq q_0$,
we have
\[
f(q) = \sum_{(q',d) \in U_q} t^d f(q')^p,
\]
whereas if $q = q_0$, we have
\[
f(q_0) = 1 + \sum_{(q',d) \in U_q} t^d f(q')^p.
\]
By Lemma~\ref{lem:vector}, $f(q)$ is algebraic over
$\FF_q(t)$ for each preradix state $q$.

For $x \in S_p \cap [0,1)$, let $s''(x)$ be the base $p$ expansion of
$x$ minus the initial radix point.
For each postradix state $q \in Q$, let $V_q$ be the set of $x \in S_p \cap
[0,1)$ such that $\delta^*(q, s''(x))$ is a final state,
and put $g(q) = \sum_{i \in V_q} t^i$. Then
if $q$ is non-final, we have
\[
g(q)^p = \sum_{d=0}^{p-1} t^d g(\delta(q,d)),
\]
whereas if $q$ is final, then
\[
g(q)^p = 1 + \sum_{d=0}^{p-1} t^d g(\delta(q,d)).
\]
By Lemma~\ref{lem:vector}, $g(q)$ is algebraic for each postradix
state $q$.

Finally, note that
\[
\sum_{i \in S} t^i = 
\sum_{q,q'} f(q) g(q'),
\]
the sum running over preradix $q$ and postradix $q'$. This
sum is algebraic over $\FF_q(t)$ by Lemma~\ref{lem:subring},
as desired.
\end{proof}
\begin{prop} \label{P:auttoalg}
Let $\sum_i x_i t^i \in \FF_q((t^\QQ))$ be a $p$-quasi-automatic
generalized Laurent series. Then $\sum_i x_i t^i$ is algebraic
over $\FF_q(t)$.
\end{prop}
\begin{proof}
Choose integers $a,b$ as in Definition~\ref{D:automatic series}.
For each $\alpha \in \FF_q$, let $S_\alpha$ be the set of
$j \in \QQ$ such that $x_{(j-b)/a} = \alpha$. Then each $S_\alpha$
is $p$-regular, so Lemma~\ref{L:auttoalg} implies that 
$\sum_{j \in S_\alpha} t^j$ is algebraic over $\FF_q(t)$.
By Lemma~\ref{lem:subring}, 
\[
\sum_i x_{ai+b} t^i = \sum_{\alpha \in \FF_q} \alpha
\left( \sum_{j \in S_\alpha} t^j \right)
\]
is also algebraic over $\FF_q(t)$; by Lemma~\ref{L:decimate},
$\sum_i x_i t^i$ is also algebraic over $\FF_q(t)$.
\end{proof}

\subsection{Algebraic implies automatic}
\label{subsec:algtoaut1}

We next prove the ``algebraic implies automatic'' direction of
Theorem~\ref{thm:main}. Unfortunately, the techniques originally
used to prove Christol's theorem (as in \cite[Chapter~12]{as})
do not suffice to give a proof of this direction. In this section,
we will get around this by using the characterization of the algebraic
closure of $\FF_q((t))$ within $\FF_q((t^\QQ))$ provided by
\cite{me}. This proof thus inherits the property of \cite{me}
of being a bit abstract, as \cite{me} uses some Galois theory and properties of
finite extensions of fields in positive characteristic (namely
Artin-Schreier theory, which comes from an argument in Galois
cohomology). It also requires invoking the ``algebraic implies
automatic'' direction of Christol's theorem itself.
We will give a second, more computationally explicit proof
of this direction later (Proposition~\ref{P:algtoaut2}).

\begin{defn}
For $c$ a nonnegative integer, let $T_c$ be the subset of $S_p$ given by
\[
T_{c} = \left\{ n - b_1 p^{-1} - b_2 p^{-2} - \cdots:
n \in \ZZ_{\geq 0}, b_i \in \{0, \dots, p-1\}, \sum b_i \leq c \right\}.
\]
\end{defn}
Then \cite[Theorem~15]{me} gives a criterion for algebraicity
of a generalized power series not over the rational function field
$\FF_q(t)$, but over the Laurent series field $\FF_q((t))$. 
It can be stated as follows.
\begin{prop} \label{prop:criterion}
For $x = \sum_i x_i t^i \in \FF_q((t^\QQ))$,
$x$ is algebraic over $\FF_q((t))$
if and only if the following conditions hold.
\begin{enumerate}
\item[(a)]
There exist integers $a,b,c \geq 0$ such that the support of
$\sum_i x_{(i-b)/a} t^i$ is contained in $T_c$.
\item[(b)]
For some $a,b,c$ as in (a),
there exist positive integers $M$ and $N$ such that
every sequence $\{c_n\}_{n=0}^\infty$ of the form
\begin{equation} \label{eq:index}
c_n = x_{(m-b- b_1 p^{-1} - \cdots - b_{j-1} p^{-j+1} -
p^{-n}(b_j p^{-j} + \cdots))/a },
\end{equation}
with $j$ a nonnegative integer, $m$ a positive integer,
and $b_i \in \{0, \dots, p-1\}$ such that
$\sum b_i \leq c$, becomes eventually periodic with period length dividing
$N$ after at most $M$ terms.
\end{enumerate} 
Moreover, in this case, (b) holds for any $a,b,c$ as in (a). 
\end{prop}
Beware that it is possible to choose $a,b$ so that the support
of $\sum x_{(i-b)/a} t^i$ 
is contained in $S_p$ and yet not have (a) satisfied
for any choice of $c$. For example,
the support of $x = \sum_{i=0}^\infty t^{(1 - p^{-i})/(p-1)}$ is contained
in $S_p$ and in $\frac{1}{p-1} T_1$, but is not contained
in $T_c$ for any $c$.

We first treat a special case of the ``algebraic implies automatic''
implication which is orthogonal to Christol's theorem.
\begin{lemma} \label{L:truncation}
Suppose that $x = \sum_i x_i t^i \in \FF_q((t^\QQ))$ has support
in $(0,1] \cap T_c$ for some nonnegative integer $c$,
and that $x$ is algebraic over $\FF_q((t))$. Then:
\begin{enumerate}
\item[(a)] $x$ is $p$-automatic;
\item[(b)] $x$ is algebraic over $\FF_q(t)$;
\item[(c)] $x$ lies in a finite set determined by $q$ and $c$.
\end{enumerate}
\end{lemma}
\begin{proof}
Note that (b) follows from (a) by virtue of
Proposition~\ref{P:auttoalg}, so it suffices to prove (a) and (c).
The criterion from Proposition~\ref{prop:criterion} applies
with $a=1,b=0$ and the given value of $c$,
so we have that every sequence $\{c_n\}_{n=0}^\infty$
of the form
\begin{equation} \label{eq:index2}
c_n = x_{1- b_1 p^{-1} - \cdots - b_{j-1} p^{-j+1} -
p^{-n}(b_j p^{-j} + \cdots) },
\end{equation}
with $b_i \in \{0, \dots, p-1\}$ such that
$\sum b_i \leq c$, becomes eventually periodic with period length dividing
$N$ after at most $M$ terms.

Define the equivalence relation on $T_c$ as follows. Declare two
elements of $T_c$ are equivalent if one can obtain the base $b$
expansion of one from the base $b$ expansion of the other by
repeating the following operation: replace a consecutive string of
$M+u+vN$ zeroes by a consecutive string of $M+u+wN$ zeroes,
where $u,v,w$ may be any nonnegative integers. The criterion
of Proposition~\ref{prop:criterion} then asserts that if $i,j
\in T_c$ satisfy $i \sim j$, then $x_{1-i} \sim x_{1-j}$; also,
the equivalence relation is clearly stable under concatenation
with a fixed postscript.

Under this equivalence relation, each equivalence class has
a unique shortest element, namely the one in which no
nonzero digit in the base $b$ expansion is preceded by $M+N$
zeroes. On one hand, this means that $x$ is determined by
finitely many coefficients, so (c) follows.
On the other hand, 
by the Myhill-Nerode theorem (Lemma~\ref{lem:mn}), it follows
that the function $f: S_p \to \FF_q$ given by $f(i) = x_{1-i}$
is $p$-automatic. (More precisely, Lemma~\ref{lem:mn} implies
that the inverse image of each element of $\FF_q$ under $f$
is $p$-regular, and hence $f$ is $p$-automatic.)
Since there is an obvious transducer that perfoms the operation
$i \mapsto 1-i$ on the valid base $b$ expansions of elements of
$S_p \cap (0,1]$ (namely, transcribe the
usual hand computation), $x$ is $p$-automatic, and (a) follows.
\end{proof}

\begin{lemma} \label{L:linearly dependent}
Suppose that $x_1, \dots, x_m \in \FF_q((t^\QQ))$ all satisfy the
hypothesis of Lemma~\ref{L:truncation} for the same value of $c$,
and that $x_1, \dots, x_m$ are linearly dependent over 
$\FF_q((t))$. Then $x_1, \dots, x_m$ are also linearly dependent
over $\FF_q$.
\end{lemma}
\begin{proof}
If $x_1, \dots, x_m$ are linearly dependent over $\FF_q((t))$,
then by clearing denominators, we can find a nonzero linear relation
among them of the form $c_1 x_1 + \cdots + c_m x_m = 0$,
where each $c_i$ is in $\FF_q \llbracket t \rrbracket$. Write
$c_i = \sum_{j=0}^\infty c_{i,j} t^j$ for $c_{i,j} \in \FF_q$; we then have
\[
0 = \sum_{j=0}^\infty \left( \sum_{i=1}^m c_{i,j} x_i \right) t^j
\]
in $\FF_q((t^\QQ))$. However, the support of the quantity in
parentheses is contained in $(j,j+1]$; in particular, these
supports are disjoint for different $j$. Thus for the sum to be
zero, the summand must be zero for each $j$; that is,
$\sum_{i=1}^m c_{i,j} x_i = 0$ for each $j$. The $c_{i,j}$ cannot
all be zero or else $c_1, \dots, c_m$ would have all been zero,
so we obtain a nontrivial linear relation among $x_1, \dots, x_m$
over $\FF_q$, as desired.
\end{proof}

We now establish the ``algebraic implies automatic'' implication
of Theorem~\ref{thm:main}.
\begin{prop} \label{P:algtoaut1}
Let $x = \sum x_i t^i \in \FF_q((t^\QQ))$ be a generalized power
series which is algebraic over $\FF_q(t)$. Then $x$ is
$p$-quasi-automatic.
\end{prop}
\begin{proof}
Choose $a,b,c$ as in Proposition~\ref{prop:criterion},
and put $y_i = x_{(i-b)/a}$ and $y = \sum_i y_i t^i$, so that
$y$ is algebraic over $\FF_q(t)$ (by Lemma~\ref{L:decimate})
and has support in $T_c$. Note that for any positive integer $m$,
$y^{q^m}$ is also algebraic over $\FF_q(t)$ and also has support
in $T_c$.
By Lemma~\ref{L:ore}, we can find a polynomial
$P(z) = \sum_{i=0}^m c_i z^{q^i}$ over $\FF_q(t)$ such that $P(y-y_0) = 0$.
We may assume without loss of generality that $c_m \neq 0$, and that
$c_l = 1$, where $l$ is the
smallest nonnegative integer for which $c_l \neq 0$. 

Let $V$ be the set of elements of $\FF_q((t^\QQ))$ which satisfy
the hypotheses of Lemma~\ref{L:truncation}; then $V$ is a finite
set which is a vector space over $\FF_q$, each of whose elements
is $p$-automatic and also algebraic over $\FF_q(t)$. Let
$v_1, \dots, v_r$ be a basis of $V$ over $\FF_q$; 
by Lemma~\ref{L:linearly dependent}, $v_1, \dots, v_r$ are also
linearly independent over $\FF_q((t))$.

By the criterion of Proposition~\ref{prop:criterion}, we can
write $y-y_0 = \sum_{j=0}^\infty v_j t^j$ with each $v_j \in
V$. In other words, $y$ is an $\FF_q((t))$-linear combination
of elements of $V$. Likewise, for each positive integer $m$,
$(y-y_0)^{q^m}$ is an $\FF_q((t))$-linear combination of elements of $V$.
For $i=l, \dots, m$, write $y^{q^i} = \sum_{j=1}^r a_{i,j} v_j$ with 
$a_{i,j} \in \FF_q((t))$; then the $a_{i,j}$ are uniquely determined by
Lemma~\ref{L:linearly dependent}.

By the same reasoning, for $j=1, \dots, r$, we can write
\[
v_h^q = \sum_{h=1}^r b_{h,j} v_j
\]
for some $b_{h,j} \in \FF_q[t]$. This means that
\begin{equation} \label{eq:system1}
a_{i,j} = \sum_h b_{h,j} a_{i-1,h}^q \qquad (i=l+1, \dots, m; j = 1, \dots, r).
\end{equation}
Moreover, the equation
\[
(y-y_0)^{q^l} = - c_{l+1} ((y-y_0)^{q^l})^q - \cdots - c_m 
((y-y_0)^{q^{m-1}})^q = 0,
\]
which holds because $P(y-y_0) = 0$ 
and $c_l = 1$ by hypothesis, can be rewritten
as
\begin{align*}
\sum_{j=1}^r a_{l,j} v_j &= \sum_{i=l}^{m-1} -c_{i+1} 
(\sum_{h=1}^r a_{i,h} v_h)^q \\
&= \sum_{i=l}^{m-1} -c_{i+1} \sum_{h=1}^r \sum_{j=1}^r
a_{i,h}^q b_{h,j} v_j.
\end{align*}
Equating coefficients of $v_j$ yields a system of equations of the form
\begin{equation} \label{eq:system2}
a_{l,j} = \sum_{g=l}^{m-1} \sum_{h=1}^r d_{g,h} a_{g,h}^q
\end{equation}
with each $d_{g,h} \in \FF_q(t)$.

Combining \eqref{eq:system1} and~\eqref{eq:system2} yields a system of
equations which translates into a matrix equation 
of the form described by Lemma~\ref{lem:vector}, in which
the matrix $B$ described therein is the identity. By Lemma~\ref{lem:vector},
the $a_{i,j}$ are algebraic over $\FF_q(t)$ for $i=l, \dots, m$.

Since each $a_{i,j}$ belongs to $\FF_q((t))$, Christol's theorem
(Theorem~\ref{T:christol}) implies that $a_{i,j}$ is $p$-automatic
for $i=l, \dots, m-1$ and $j=1, \dots, r$. This implies that
$(y-y_0)^{q^l} = \sum_{j=1}^r a_{l,j} v_j$ is
$p$-automatic, as follows. Write 
$(y-y_0)^{q^l} = \sum_{k=0}^\infty w_k t^k$ with
$w_k \in V$. Then the function $k \mapsto w_k$ is $p$-automatic
because the $a_{i,j}$ are $p$-automatic. So we can build an automaton
that, given a base $b$ expansion, sorts the preradix string
$k$ according to the value of $w_k$, then handles the postradix string
by imitating some automaton that computes $w_k$.

Since $(y-y_0)^{q^l}$ is $p$-automatic, $y-y_0$ is $p$-quasi-automatic,
as then is $y$, as then is $x$.
This yields the desired result.
\end{proof}

Note that Propositions~\ref{P:auttoalg} and~\ref{P:algtoaut1} together
give a complete proof of Theorem~\ref{thm:main}. We will give a
second proof of a statement equivalent to Proposition~\ref{P:algtoaut1}
later (see Proposition~\ref{P:algtoaut2}).

\section{Polynomials over valued fields}
\label{sec:valued}

In this chapter, we introduce some additional algebraic tools that
will help us give a more algorithmic proof of the
``algebraic implies automatic'' implication of Theorem~\ref{thm:main}.
Nothing in this chapter is particularly novel, but the material may
not be familiar to nonspecialists, so we give a detailed presentation.

\subsection{Twisted polynomial rings}
\label{subsec:twisted}

Ore's lemma (Lemma \ref{L:ore}) asserts that every
algebraic element of a field of characteristic $p$ over a subfield
is a root of an additive polynomial. It is sometimes more convenient
to view additive polynomials as the result of applying ``twisted
polynomials'' in the Frobenius operator. These polynomials arise naturally
in the theory of Drinfeld modules; see for instance \cite{hayes}.

Throughout this section, the field $K$ will have characteristic $p>0$.

\begin{defn}
Let $K \{F\}$ denote
the noncommutative ring whose elements are finite formal sums
$\sum_{i=0}^m c_i F^i$, added componentwise and multiplied by the
rule
\[
\left( \sum_{i=0}^m c_i F^i \right)
\left( \sum_{j=0}^n d_j F^j \right)
= \sum_{k=0}^{m+n} \left( \sum_{i+j=k} c_i d_j^{p^i} \right) F^k.
\]
The ring $K \{F\}$ is called the \emph{twisted polynomial ring} over $K$.
(Note that the same definition can be made replacing the $p$-power 
endomorphism by any endomorphism of $K$, but we will only use this 
particular form of the construction.)
\end{defn}

\begin{defn}
As in the polynomial case, the \emph{degree} of a nonzero twisted polynomial
$\sum c_i F^i$ is the largest $i$ such that $c_i \neq 0$; we conventionally
take the degree of the zero polynomial to be $-\infty$.
The degree of the product of two nonzero twisted polynomials is the
sum of their individual degrees; in particular, the ring $K\{F\}$ is an
integral domain.
\end{defn}

Twisted polynomial rings admit the following ``right division algorithm'',
just as in the usual polynomial case.
\begin{lemma}
Let $S(F)$ and $T(F)$ be twisted polynomials over $K$, with $T(F)$ nonzero
and $\deg(T) = d \geq 0$. Then there exists a unique pair $Q(F), R(F)$
of twisted polynomials over $K$ with $\deg(R) < d$, such that
$S = QT + R$.
\end{lemma}
\begin{proof}
Write $T(F) = \sum_{i=0}^d c_i F^i$ with $c_d \neq 0$.
Existence follows by induction on $\deg(S)$ and the fact that we can construct
a left multiple of $T$ of any prescribed degree $e \geq d$ with any
prescribed leading coefficient $a$ (namely $(a/c_d^{p^{e-d}}) F^{e-d}$).
Uniqueness follows from the fact that if $S = QT + R = Q'T + R'$ are
two decompositions of the desired form, then $R-R' = (Q'-Q)T$ is a left 
multiple of $T$ but $\deg(R-R') < \deg(T)$, so $R-R' = 0$.
\end{proof}

We may view twisted polynomials over $K$ as additive operators on any field
$L$ containing $K$ by declaring that
\[
\left( \sum_i c_i F^i \right)(z) = \sum_i c_i z^{p^i}.
\]
\begin{lemma} \label{L:Fp-vector space}
Let $L$ be an algebraic closure of $K$, and let
$T(F) = \sum_{i=0}^d c_i F^i$ be a nonzero twisted polynomial of degree $d$
over $K$.
Then the kernel $\ker_L(T)$ of $T$ acting on $L$ is an $\FF_p$-vector
space of dimension $\leq d$, with equality if and only if $c_0 \neq 0$.
\end{lemma}
\begin{proof}
The kernel is an $\FF_p$-vector space because $T(F)$ is an additive operator
on $L$ (Lemma~\ref{L:additive}), 
and the dimension bound holds because $T(F)(z)$ is a polynomial
in $z$ of degree $p^d$. The equality case holds because
the formal derivative in $z$ of $T(F)(z)$ is $c_0$, and this derivative
vanishes if and only if $T(F)(z)$ has no repeated roots over $L$.
\end{proof}

\begin{lemma} \label{L:left multiple}
Let $S(F)$ and $T(F)$ be twisted polynomials over $K$, with the constant
coefficient of $T$ nonzero, and let $L$ be an
algebraic closure of $K$. Then
$S$ is a left multiple of $T$ (that is, $S = QT$ for some $Q \in K\{F\}$)
if and only if $\ker_L(T) \subseteq \ker_L(S)$.
\end{lemma}
\begin{proof}
If $S = QT$, then $T(F)(z) = 0$ implies $S(F)(z) = 0$, so
$\ker_L(T) \subseteq \ker_L(S)$. Conversely, suppose that
$\ker_L(T) \subseteq \ker_L(S)$, and write $S = QT + R$ by the right
division algorithm; we then have $\ker_L(T) \subseteq \ker_L(R)$ as well.
By Lemma~\ref{L:Fp-vector space}, $\ker_L(T)$ is an $\FF_p$-vector
space of dimension $\deg(T)$; if $R$ is nonzero, then 
$\ker_L(R)$ is an $\FF_p$-vector space of dimension
at most $\deg(R) < \deg(T)$. But that would contradict the inclusion
$\ker_L(T) \subseteq \ker_L(R)$, so we must have $R=0$ and $S = QT$,
as desired.
\end{proof}

\subsection{Newton polygons}
\label{subsec:np}

The theory of Newton polygons is a critical ingredient in the computational
and theoretical study of valued fields. We recall a bit of this theory here.

\begin{defn}
For $x \in \FF_q((t^\QQ))$ nonzero, let $v(x)$ denote the smallest
element of the support of $x$; we call $v$ the \emph{valuation} on
$\FF_q((t^\QQ))$. We also formally put $v(0) = \infty$. The function
$v$ has the usual properties of a valuation:
\begin{align*}
v(x+y) &\geq \min\{v(x), v(y)\} \\
v(xy) &= v(x) + v(y).
\end{align*}
Given a nonzero polynomial $P(z) = \sum_i c_i z^i$
over $\FF_q((t^\QQ))$,
we define the \emph{Newton polygon} of $P$ to be the lower boundary
of the lower convex hull of the set of points $(-i, v(c_i))$.
The slopes of this polygon are called the \emph{slopes} of $P$;
for $r \in \QQ$, we define the \emph{multiplicity} of $r$ as a
slope of $P$ to be the width (difference in $x$-coordinates between
the endpoints) of the segment of the Newton polygon of $P$ of slope
$r$, or 0 if no such segment exists.
We say $P$ is \emph{pure (of slope $r$)} if all of the slopes of $P$
are equal to $r$.
We conventionally declare that $\infty$ may also be a slope,
and its multiplicity as a slope of a polynomial $P$ is the order of
vanishing of $P$ at $z=0$.
\end{defn}

\begin{lemma} \label{L:slopes add}
Let $P$ and $Q$ be nonzero polynomials over $\FF_q((t^\QQ))$. 
Then for each $r \in \QQ \cup \{\infty\}$, 
the multiplicity of $r$ as a slope of $P+Q$
is the sum of the multiplicities of $r$ as a slope of $P$ and of $Q$.
\end{lemma}
\begin{proof}
The case $r = \infty$ is clear, so we assume $r \in \QQ$.
Write $P(z) = \sum_i c_i z^i$ and $Q(z) = \sum_j d_j z^j$.
Let $(-e,v(c_e))$ and $(-f, v(d_f))$ (resp.\
$(-g,v(d_g))$ and $(-h,v(d_h))$) be the left and right endpoints,
respectively, of the (possibly degenerate) segment in which the
Newton polygon of $P$ (resp.\ of $Q$) meets its support line of slope $r$.
Then 
\begin{gather*}
v(c_i) + ri \geq v(c_e) + re = v(c_f) + rf \\
v(d_j) + rj \geq v(d_g) + rg = v(d_h) + rh
\end{gather*}
with strict inequality if $i \notin [e,f]$ or $j \notin [g,h]$. For
each $k$,
\[
\min_{i+j=k} \{ v(c_i d_j) + rk \} \geq v(c_e) + re + v(d_g) + rg
\]
and the inequality is strict in each of the following cases:
\begin{itemize}
\item
$k \notin [f+h,e+g]$;
\item
$k = e+g$ and $i \neq e$;
\item
$k = f+h$ and $i \neq h$.
\end{itemize}
If we write $R = PQ$ and $R(z) = \sum a_k z^k$, it follows that
\[
\min_k \{ v(a_k) + rk \} \geq v(c_e) + re + v(d_g) + rg,
\]
with equality for $k = f+h$ and $k = e+g$ but not for any
$k \notin [f+h,e+g]$. It follows that the multiplicity of $r$ as 
a slope of $R$ is $e+g-(f+h) = (e-f) + (g-h)$, proving the desired
result. 
\end{proof}

\begin{cor}
Let $P(z)$ be a nonzero polynomial
over $\FF_q((t^\QQ))$ which factors as
a product $Q_1 \cdots Q_n$ of pure polynomials (e.g., linear polynomials). 
Then for each $r \in \QQ \cup \{\infty\}$,
the sum of the degrees of the $Q_i$, over those $i$ for which $Q_i$ has
slope $r$, is equal to the multiplicity of $r$ as 
a slope of $P$.
\end{cor}

\subsection{Slope splittings}
\label{subsec:slope-split}

We now recall a special form of Hensel's lemma, which makes it possible
to split polynomials by their slopes.

\begin{defn}
A subfield $K$ of $\FF_q((t^\QQ))$ is \emph{closed} under the valuation
$v$ if for any sequence $\{z_n\}_{n=0}^\infty$ such that
$v(z_n - z_{n+1}) \to \infty$ as $n \to \infty$, there exists $z \in K$
such that $v(z_n - z) \to \infty$ as $n \to \infty$. For instance,
$\FF_q((t))$ is closed, as is any finite extension of $\FF_q((t))$.
\end{defn}

\begin{lemma} \label{L:splitting}
Let $K$ be a subfield of $\FF_q((t^\QQ))$ closed under the valuation $v$.
Let $P(z)$ be a nonzero polynomial over $K$, and choose $u \in \QQ$ which does
not occur as a slope of $P$.
Then there exists a factorization $P = QR$, where $Q$ is a polynomial
over $K$ with all slopes less than $u$, and $R$ is a polynomial over
$K$ with all slopes greater than $u$.
\end{lemma}
\begin{proof}
Pick $r,s \in \QQ$ such that $r<s$ and the interval $[r,s]$ contains $u$ but
does not contain any slopes of $P$.
Let $S$ be the ring of formal sums
$\sum_{j \in \ZZ} a_j z^j$ over $K$
such that $v(a_j) + rj \geq 0$ for $j \geq 0$
and $v(a_j) + sj \geq 0$ for $j \leq 0$.
Define the valuation $v_u$ on $S$ by
the formula
\[
v_u\left(\sum a_j z^j \right) = \min_j \{v(a_j) + uj\}.
\]
Given $x = \sum a_j z^j \in S$, write $f(x) = \sum_{j > 0} a_j z^j$,
so that $f(f(x)) = f(x)$.

Write $P(z) = \sum_i c_i z^i$, put $d = \deg P$,
let $e$ be the unique integer
that minimizes $v(c_e) + ue$, and put
$x_0 = c_e^{-1} z^{-e} P(z) \in S$.
Define the sequence $\{x_h\}_{h=0}^\infty$ by the recurrence
$x_{h+1} = x_h (1-f(x_h))$ for $h = 0,1,\dots$.

Put $\ell = v_u(x_0 - 1)$, so that $\ell > 0$. We now show by induction
that $v_u(x_h - 1) \geq \ell$ and $v_u(f(x_h)) \geq (h+1)\ell$
for all $h$. The first inequality clearly implies the second for $h=0$;
given both inequalities for some $h$, we have
\begin{align*}
v_u(x_{h+1} - 1) &= v_u(x_h - x_h f(x_h) - 1) \\
&\geq \min\{v_u(x_h - 1), v_u(x_h) + v_u(f(x_h))\} \\
&\geq \min\{\ell, (h+1)\ell\} \\
&\geq \ell
\end{align*}
and
\begin{align*}
v_u(f(x_{h+1})) &= v_u(f(x_h - x_h f(x_h))) \\
&= v_u(f(x_h - f(x_h) - (x_h-1)f(x_h))) \\
&= v_u(f(x_h) - f(f(x_h)) - f((x_h-1)f(x_h))) \\
&= v_u(f((x_h-1)f(x_h))) \\
&\geq v_u((x_h-1)f(x_h)) \\
& \geq \ell + (h+1) \ell \\
&= (h+2) \ell.
\end{align*}
Thus the sequence $\{x_h\}$ converges to a limit $x \in S$.

By construction, $f(x) = 0$, so $x$ has no positive powers of $z$.
By induction on $h$,
$z^e x_h$ has no negative powers of $z$ for each $h$, so
$z^e x$ also has no negative powers of $z$. Thus 
$z^e x$ is a polynomial in $z$ of degree exactly $e$,
and its slopes are all less than or equal to $r$.

On the other hand, we have
$x = x_0(1-f(x_0))(1-f(x_1))\cdots$,
so that $x_0 x^{-1}$ has only positive powers of $z$.
However, $z^{e-d}x_0 x^{-1}$ has no positive powers of $z$, so
$x_0 x^{-1}$ is a polynomial in $z$ of degree exactly $d-e$,
and its slopes are all greater than or equal to $s$.

We thus obtain the desired factorization
as $P = QR$ with $Q(z)= z^e x$ and $R(z) = c_e x_0 x^{-1}$.
\end{proof}

\begin{cor} \label{C:slope factorization}
Let $K$ be a subfield of $\FF_q((t^\QQ))$ closed under the valuation $v$.
Then every monic polynomial $P$ over $K$ admits a unique factorization
$Q_1 \cdots Q_n$ into monic polynomials, such that each $Q_i$ is pure of 
some slope $s_i$, and $s_1 < s_2 < \cdots < s_n$.
\end{cor}
\begin{proof}
The existence of the factorization follows from Lemma~\ref{L:splitting};
the uniqueness follows from unique factorization for polynomials over $K$
and the additivity of multiplicities (Lemma~\ref{L:slopes add}).
\end{proof}
\begin{defn}
We call the factorization given by Corollary~\ref{C:slope factorization}
the \emph{slope factorization} of the polynomial $P$.
\end{defn}

\subsection{Slope splittings for twisted polynomials}
\label{subsec:slope-twisted}

\begin{remark}
Throughout this section, we will use the fact that for $K$ any subfield of
$\FF_q((t^\QQ))$, the valuation $v$ extends uniquely to any algebraic
closure of $K$ \cite[Proposition~II.3]{ser}
We do this for conceptual clarity; not doing so
(at the expense of making the arguments messier and more computational)
would not shed any light on how to compute with automatic series.
\end{remark}

\begin{lemma} \label{L:also additive}
Let $K$ be a subfield of $\FF_q((t^\QQ))$ closed under the valuation $v$.
Let $P$ be a monic additive polynomial over $K$,
and let $Q_1 \cdots Q_n t^{p^d}$ be its slope factorization.
Then $Q_n t^{p^d}$ is also additive.
\end{lemma}
\begin{proof}
Let $L$ be an algebraic closure of $K$.
Let $s_1, \dots, s_n$ be the slopes of $Q_n$,
and let $V$ be the roots of $P$ in $L$; then the possible valuations of the
nonzero elements of $V$ are precisely $s_1, \dots, s_n$.
Moreover, an element of $V$ is a root of $Q_n$ if and only if
$v(x) \geq s_n$, and this subset of $V$ is an $\FF_p$-subspace of $V$.
By Lemma~\ref{L:additive}, $Q_n t^{p^d}$ is also additive.
\end{proof}

In terms of twisted polynomials, we obtain the following analogue
of the slope factorization.
\begin{defn}
Let $P(F)$ be a twisted polynomial over $\FF_q((t^\QQ))$.
For $r \in \QQ$, we say $P(F)$ is \emph{pure of slope $r$}
if $P$ has nonzero constant term and
the ordinary polynomial $P(F)(z)/z$ is pure of slope $r$.
We conventionally say that any power of $F$ is pure of slope $\infty$.
\end{defn}

\begin{prop} \label{P:split twisted}
Let $K$ be a subfield of $\FF_q((t^\QQ))$ closed under the valuation $v$,
and let $P(F)$ be a monic twisted polynomial over $K$.
Then there exists a factorization $P = Q_1 \cdots Q_n$ of $P$
into monic twisted polynomials over $K$, in which each $Q_i$ is pure of some
slope.
\end{prop}
\begin{proof}
Let $R(z)$ be the highest finite slope factor in the slope factorization
of $P(F)(z)$, times the slope $\infty$ factor (a power of $z$). 
By Lemma~\ref{L:also additive}, 
$R(z)$ is additive, so we have $R(z) = Q(F)(z)$ for some twisted
polynomial $Q$. By Lemma~\ref{L:left multiple}, we can factor
$P = P_1 Q$ for some $P_1$ of lower degree than $P$; repeating the
argument yields the claim.
\end{proof}

We can split further at the expense of enlarging $q$.
\begin{prop} \label{P:split twisted2}
Let $K$ be a
subfield of $\FF_q((t^\QQ))$ closed under the valuation $v$,
and let $P(F)$ be a monic twisted polynomial over $K$ which is pure
of slope $0$. Then the polynomial $P(F)(z)$ factors completely (into
linear factors) over
$K \otimes_{\FF_q} \FF_{q'}$, for some power $q'$ of $q$.
\end{prop}
\begin{proof}
Write $P(F) = \sum_{i=0}^d c_i F^i$ with $c_d = 1$. Since
$P$ is pure of slope 0, we have $v(c_0) = 0$
and $v(c_i) \geq 0$ for $1 \leq i \leq d-1$.

Let $a_i \in \FF_q$ be the constant coefficient of $c_i$; then for
some power $q'$ of $q$, the polynomial $\sum a_i z^{p^i}$ has $p^d$
distinct roots in $\FF_{q'}$. Let $r \in \FF_{q'}$ be a nonzero
root of $\sum a_i z^{p^i}$; then $P(F)(z+r)$ has one slope greater than
$0$ and all others equal to $0$. The slope factorization of $P(F)(z+r)$
then has a linear factor; in other words, $P(F)(z)$ has a unique
root $r' \in K \otimes_{\FF_q} \FF_{q'}$ with $v(r-r') > 0$.
By the same argument, $P(F)(z)$ factors completely over
$K \otimes_{\FF_q} \FF_{q'}$, as desired.
\end{proof}
\begin{cor} \label{C:splits completely}
Let $K$ be a subfield of $\FF_q((t^\QQ))$, containing all fractional
powers of $t$ and
closed under the valuation $v$,
and let $P(F)$ be a monic twisted polynomial over $K$.
Then for some power $q'$ of $q$, there exists a factorization
$P = Q_1 \cdots Q_n$ of $P$ into monic \emph{linear} twisted polynomials
over $K \otimes_{\FF_q} \FF_{q'}$.
\end{cor}
\begin{proof}
We may proceed by induction on $\deg P$; it suffices to show that
if $P$ is not linear, then it is a left multiple of some linear
twisted polynomial over $K \otimes_{\FF_q} \FF_{q'}$ for some
$q'$.
By Proposition~\ref{P:split twisted}, we may reduce to the case
where $P$ is pure of some slope; by rescaling the polynomial $P(F)(z)$
(and using the fact that $K$ contains all fractional powers of $t$),
we may reduce to the case where $P$ is pure of slope 0.
By Proposition~\ref{P:split twisted2}, $P(F)(z)$ then splits completely 
over $K \otimes_{\FF_q} \FF_{q'}$ for some $q'$. Choose any one-dimensional
$\FF_p$-subspace of the set of roots of $P(F)(z)$; by
Lemma~\ref{L:additive}, these form
the roots of $Q(F)(z)$ for some monic linear twisted polynomial $Q$
over $K \otimes_{\FF_q} \FF_{q'}$.
By Lemma~\ref{L:left multiple}, we can write $P = P_0 Q$ for some
$P_0$; this completes the induction and yields the desired result.
\end{proof}

\section{Proof of the main theorem: concrete approach}
\label{sec:concrete}

In this chapter, we study the properties of automatic power series more
closely, with the end of giving a more down-to-earth proof (free of
any dependence on Galois theory or the like) of the
``algebraic implies automatic'' implication of Theorem~\ref{thm:main}.
In so doing, we introduce some techniques which may be of use in
explicitly computing in the algebraic closure of $\FF_q(t)$; however,
we have not made any attempt to ``practicalize'' these techniques.
Just how efficiently one can do such computing is a problem worthy
of further study; we discuss this question briefly in the next chapter.

\subsection{Transition digraphs of automata}
\label{subsec:transition}

We begin the chapter with 
a more concrete study of the automata that give rise to
generalized power series.

\begin{defn}
Let $M$ be a DFAO.
Given a state $q_1 \in Q$,
a state $q \in Q$ is said to be \emph{reachable from $q_1$} if
$\delta^*(q_1,s) = q$ for some string $s \in \Sigma^*$, and
\emph{unreachable from $q_1$} otherwise; if $q_1 = q_0$, we
simply say that $q$ is reachable or unreachable.
Any state from which a state unreachable from $q_0$ can be reached
must itself be unreachable from $q_0$.
A DFAO with no unreachable states
is said to be \emph{minimal}; given any DFAO, one can remove all of its
unreachable 
states 
to obtain a minimal DFAO that accepts the same language.
\end{defn}

\begin{defn}
A state $q \in Q$ is said to be \emph{relevant} if there exists
a final state reachable from $q$, and \emph{irrelevant} otherwise.
Any state which is reachable from an irrelevant state
must itself be irrelevant.
\end{defn}

\begin{defn}
Let $M$ be a DFA with input alphabet $\Sigma_b = \{0, \dots, b-1, .\}$.
We say $M$ is \emph{well-formed} (resp.\ \emph{well-ordered})
if the language accepted by $M$
consists of the valid base $b$ expansions
of the elements of an arbitrary (resp.\ a well-ordered) 
subset of $S_b$.
\end{defn}

The goal of this section is to characterize minimal well-ordered DFAs
via their transition graphs; to do so, we use the following definitions.
(By way of motivation,
a connected undirected graph is called a ``cactus'' if each vertex
of the graph lies on exactly one minimal cycle. In real life, the saguaro
is a specific type of cactus indigenous to the southwest United States and northwest Mexico.)
\begin{defn}
A directed graph $G = (V,E)$ equipped with a distinguished vertex $v \in V$
is called a \emph{rooted saguaro}
if it satisfies the following conditions.
\begin{enumerate}
\item[(a)] Each vertex of $G$ lies on at most one minimal cycle.
\item[(b)] There exist directed paths from $v$ to each vertex of $G$.
\end{enumerate}
In this case, we say that $v$ is a \emph{root} of the saguaro.
(Note that a saguaro may have more than one root, because any vertex
lying on the same minimal cycle as a root is also a root.)
A minimal cycle of a rooted saguaro is called a \emph{lobe}.
An edge of a rooted saguaro is \emph{cyclic} if it lies on a lobe
and \emph{acyclic} otherwise.
\end{defn}

\begin{defn}
Let $G = (V,E)$ 
be a rooted saguaro. A \emph{proper $b$-labeling} of $G$ is a function
$\ell: E \to \{0, \dots, b-1\}$ with the following properties.
\begin{enumerate}
\item[(a)]
If $v,w,x \in V$ and $vw, vx \in E$, then $\ell(vw) \neq \ell(vx)$.
\item[(b)]
If $v,w,x \in V$, $vw \in E$ lies on a lobe, and
$vx \in E$ does not lie on a lobe, then $\ell(vw) > \ell(vx)$.
\end{enumerate}
\end{defn}

\begin{theorem}
Let $M$ be a DFA with input alphabet $\Sigma_b$, which is 
minimal and well-formed.
Then $M$ is well-ordered if and only if for each relevant postradix state
$q$, the subgraph $G_q$
of the transition graph consisting of relevant states
reachable from $q$ is a rooted saguaro with root $q$, equipped
with a proper $b$-labeling.
\end{theorem}
\begin{proof}
First suppose that $M$ is well-ordered. Let $q$ be a 
relevant postradix state of $M$; note that all transitions from $q$ are
labeled by elements of $\{0, \dots, b-1\}$, since a valid base $b$
expansion cannot have two radix points. 
Suppose that $q$ admits a cyclic transition by $s \in \{0, \dots, b-1\}$
and also admits a transition to a relevant
state by $s' \neq s$. We will show that $s > s'$.

To see this, choose $w \in \Sigma^*$ of minimal
length such that $\delta^*(q,s w)
= q$. Since $M$ is minimal, $q$ is reachable, so
we can choose $w_0 \in \Sigma^*$ such that
$\delta^*(q_0,w_0) = q$. Since $\delta(q,s')$ is relevant,
we can choose $w_1 \in \Sigma^*$ such that
$\delta^*(q,s'w_1) \in F$.
Then all of the strings
\[
w_0 s' w_1, w_0 s w s' w_1, w_0 s w s w s' w_1, \dots
\]
are accepted by $M$; however, if $s' > s$, these form the
valid base $b$ expansions of an infinite decreasing sequence,
which would contradict the fact that $M$ is well-formed.
Hence $s > s'$.

This implies immediately that a relevant postradix state cannot
lie on more than one minimal cycle, so $G_q$ is a rooted saguaro.
Moreover, no two edges from the same vertex of $G_q$ carry the same label,
and what we just proved implies that
any cyclic edge must carry a greater label than any
other edge from the same vertex. Hence $G_q$ carries a proper
$b$-labeling.

Suppose conversely that each $G_q$ is a rooted saguaro carrying
a proper $b$-labeling, but assume by way of contradiction
that the language accepted by $M$ includes an infinite decreasing
sequence $x_1, x_2, \dots$. We say the $m$-th digit of the base $b$ 
expansion of $x_i$ (counting from the left)
is \emph{static} if $x_i, x_{i+1}, \dots$
all have the same $m$-th digit. Then the expansion of each $x_i$ begins with
an initial segment (possibly empty) of static digits, and the number
of static digits tends to infinity with $i$. There exists a unique
infinite sequence $s_1,s_2, \dots$ of elements of $\Sigma$
such that each initial segment of static digits occurring in the
expansion of some $x_i$ has the form $s_1\cdots s_m$ for some
$m$ (depending on $i$).

The number of preradix digits in the expansion of $x_i$ is at most
the corresponding number for $x_1$, so the sequence $s_1,s_2,\dots$
must include one (and only one) radix point. Define $q_m = 
\delta^*(q_0, s_1\cdots s_m)$; then each $q_m$ is relevant, and $q_m$
is postradix for $m$ sufficiently large. Moreover, only
finitely many of the postradix transitions from $q_m$ to $q_{m+1}$ can be
acyclic, as otherwise some such transition would repeat, but then would
visibly be part of a cycle. 

Choose $x_i$ in the decreasing sequence whose first $m$ digits
are all static, for $m$ large enough that $q_m$ is postradix 
and the transition from $q_n$ to $q_{n+1}$ is cyclic for all $n \geq m$
(again, this amounts to taking $i$ sufficiently large).
Let $t_1 t_2 \cdots \in \Sigma^*$ denote the base $b$ expansion of $x_i$.
Since $x_{i+1} < x_i$, there exists a smallest integer $n>m$
such that $s_n \neq t_n$, and that integer $n$ satisfies
$s_n < t_n$. However, the transition from $q_n$ to $q_{n+1}$
is cyclic, and so the inequality $s_n < t_n$ violates the definition
of a proper $b$-labeling. This contradiction means that there
cannot exist an infinite decreasing subsequence, so
$M$ is well-ordered, as desired.
\end{proof}

\subsection{Arithmetic for well-ordered automata}
\label{subsec:arithmetic}

We next verify directly that the set of automatic series is
a ring; this follows from Theorem~\ref{thm:main} and
Lemma~\ref{lem:subring}, but giving a direct proof will on one hand
lead to a second proof of Theorem~\ref{thm:main}, and on the other
hand suggest ways to implement computations on automatic series
in practice. (However, we have not made any attempt here to optimize
the efficiency of the computations; this will require further study.)

We first note that adding automatic series is easy.
\begin{lemma} \label{lem:add automatic}
Let $x,y \in \FF_q((t^\QQ))$ be $p$-automatic (resp.\ $p$-quasi-automatic). 
Then $x+y$ is $p$-automatic (resp.\ $p$-quasi-automatic).
\end{lemma}
\begin{proof}
The claim for $p$-quasi-automatic series follows from the claim
for $p$-automatic series by Lemma~\ref{L:dilation}, so we may as well
assume that $x,y$ are $p$-automatic.
Write $x = \sum x_i t^i$ and $y = \sum y_i t^i$;
then by assumption, the functions $i \mapsto x_i$ and $i \mapsto y_i$
on $S_p$ are $p$-automatic. Hence the function
$i \mapsto (x_i,y_i)$ is also $p$-automatic, as then
is the function $i \mapsto x_i + y_i$.
This yields the claim.
\end{proof}

We next tackle the stickier subject of multiplication of automatic
series.
\begin{lemma} \label{lem:mult automatic}
Let $x,y \in \FF_q((t^\QQ))$ be $p$-automatic (resp.\ $p$-quasi-automatic). 
Then $xy$ is $p$-automatic (resp.\ $p$-quasi-automatic).
\end{lemma}
Beware that this proof works primarily with \emph{reversed} base $p$
expansions; the notions of leading and trailing zeroes will be in terms
of the reversed expansions, so leading zeroes are in the \emph{least}
significant places.
\begin{proof}
Again, it suffices to treat the automatic case.
Write each of $x$ and $y$ as an $\FF_q$-linear combination of generalized
power series of the form $\sum_{i \in S} t^i$ for some $S \subset S_p$;
by Lemma~\ref{lem:add automatic}, if the product of two series
of this form is always $p$-automatic, then so is $xy$.

That is, we may assume without loss of generality that
$x = \sum_{i \in A_1} t^i$ and $y = \sum_{i \in A_2} t^i$.
Let $S$ be the subset of $\Sigma_p^* \times \Sigma_p^*$ consisting
of pairs $(w_1, w_2)$ with the following properties.
\begin{enumerate}
\item[(a)] $w_1$ and $w_2$ have the same length.
\item[(b)] $w_1$ and $w_2$ each end with 0.
\item[(c)] $w_1$ and $w_2$ each have a single radix point, and both are in
the same position.
\item[(d)] After removing leading and trailing zeroes, $w_1$ and $w_2$ become
the reversed valid base $p$ expansions of some
$i,j \in S_p$.
\item[(e)] The pair $(i,j)$ belongs to $A_1 \times A_2$.
\end{enumerate}
By (a), we may view $S$ as a language over $\Sigma_p \times \Sigma_p$;
it is straightforward to verify that this language is in fact regular.
Let $M = (Q, \Sigma, \delta, q_0, F)$ be a DFA accepting $S$,
with $\Sigma = \Sigma_p \times \Sigma_p$.

Define an NFA $M' = (Q', \Sigma', \delta', q'_0, F')$ as follows.
Put $Q' = Q \times \{0,1\}$ and $\Sigma' = \Sigma_p$.
For $(q,i) \in Q'$ and $s \in \{0, \dots, p-1\}$, we include
$(q',0)$ (resp.\ $(q',1)$) 
in $\delta((q,i),s)$ if there exists a pair $(t,u) \in \{0,
\dots, p-1\} \times \{0, \dots, p-1\}$ with $t+u+i < p$
(resp.\ $t+u+i \geq p$) and $t+u+i \equiv s \pmod{p}$
such that $\delta(q, (t,u)) = q'$;
for $(q,i) \in Q'$ and $s$ equal to the radix point, we
include $(q',i)$ in $\delta((q,i),s)$ if 
$\delta(q,(s,s)) = q'$ (and we never include $(q',1-i)$).
Put $q'_0 = (q_0,0)$ and put $F' = F \times \{0\}$.

Suppose $w$ is a string which, upon removal of leading and trailing
zeroes, becomes
the reversed valid base $p$ expansion of some $z \in S_p$.
Then the number of accepting paths of $w$ in $M'$ is equal to the
number of pairs $(w_1, w_2) \in S$ which sum to $w$ 
\emph{with its leading and trailing zeroes} under
ordinary base $p$ addition with carries. By 
Lemma~\ref{lem:quant}, taking this number modulo $p$ yields a
finite-state function.

Suppose further that $w$ begins with $m$ leading zeroes, for 
$m$ greater than the number of states of $M$. We claim that
for any pair $(w_1, w_2) \in S$
which sums to $w$, both $w_1$ and $w_2$ must begin with a leading
zero. Namely, if this were not the case, then in processing
$(w_1, w_2) \in S$ under $M$, some state must be repeated
within the first $m$ digits. 
Let $(b_1,b_2)$ be the strings that lead to the first
arrival at the repeated state, let $(w_1,w_2)$ be the strings between
the first and second arrivals, and let $(e_1,e_2)$ be the remaining strings.
Then
\[
b_ie_i, b_iw_ie_i, b_iw_iw_ie_i, \dots,
\]
represent the reversed base $p$ expansions (with possible trailing zeroes,
but no leading zeroes) of some elements $z_{0i}, z_{1i}, \dots$ of $A_i$.
Moreover, the numbers $z_{0i}, z_{1i}, \dots$ are all distinct:
if $b_i$ is nonempty, then $b_i$ begins with a nonzero digit, while if
$b_i$ is empty, then $w_i$ begins with a nonzero digit.

We next verify that $z_{j0} + z_{j1} = z$ for each $j$.
We are given this assertion this for $j=1$;
in doing that addition, the stretch during which $w_1$ and $w_2$
are added begins with an incoming carry and ends with an outgoing carry.
Moreover, all digits produced before and during the stretch are zeroes.
Thus we may remove this stretch, or repeat it at will, without changing
the base $p$ number represented by the sum (though the number of leading
zeroes will change).

Since $z_{j0} + z_{j1} = z$ for all $j$, one of the sequences
$z_{0i}, z_{1i}, \dots$ must be strictly decreasing. This yields
a contradiction, implying that both $w_1$ and $w_2$ had to begin with
leading zeroes. We conclude from this that the number
of accepting paths of $w$ in $M'$ is preserved by adding a leading zero
to $w$ provided that $w$ already has $m$ leading zeroes.

It follows that the function that, given the reversed valid base $p$
expansion of a number $k \in S_p$, computes the mod $p$ reduction
of the number of ways to write $k=i+j$ with $i \in A_1$ and $j \in A_2$,
is a finite-state function: we may compute it by appending 
$m$ leading zeroes and one trailing zero to the reversed expansion and then
running the result through $M'$. Thus $xy$ is $p$-automatic, as desired.
\end{proof}

Division seems even more complicated to handle directly (though we suspect
it is possible to do so); we treat it here
with an indirect approach (via the ``automatic implies algebraic''
direction of Theorem~\ref{thm:main}).
\begin{lemma} \label{lem:rattoaut}
If $x \in \FF_q(t)$, then $x$ (viewed as an element of $\FF_q((t^\QQ))$)
is $p$-automatic.
\end{lemma}
\begin{proof}
There is no loss of generality in assuming $x \in \FF_q \llbracket t
\rrbracket$. Writing $x = \sum_{i=0}^\infty c_i t^i$, we see that the
sequence $\{c_i\}$ is linear recurrent over $\FF_q$, hence eventually
periodic. By \cite[Theorem~5.4.2]{as}, $x$ is $p$-automatic.
\end{proof}

\begin{prop}
The set of $p$-quasi-automatic series is a subfield of $\FF_q((t^\QQ))$
contained in the integral closure of $\FF_q(t)$.
\end{prop}
\begin{proof}
Let $S$ be the set of $p$-quasi-automatic series.
By Lemmas~\ref{lem:add automatic}
and~\ref{lem:mult automatic}, $S$ is a subring of $\FF_q((t^\QQ))$.
By Lemma~\ref{lem:rattoaut}, $S$ contains $\FF_q(t)$;
by Proposition~\ref{P:algtoaut1}, each element of $S$
is algebraic over $\FF_q(t)$. 
Hence given $x \in S$, there exists a polynomial $P(z) =
c_0 + c_1 z + \cdots + c_n z^n$ over $\FF_q(t)$ with $P(x) = 0$,
and we may assume without loss of generality that $c_0, c_n \neq 0$.
We can then write
\[
z^{-1} = c_0^{-1} (c_1 + c_2 z + \cdots + c_n z^{n-1}),
\]
and the right side is contained in $S$. We conclude that $S$ is closed under
taking reciprocals, and hence is a subfield of 
$\FF_q((t^\QQ))$ contained in the integral closure of
$\FF_q(t)$, as desired.
\end{proof}

\subsection{Newton's algorithm}
\label{subsec:Newton}

To complete the ``concrete'' proof of the ``algebraic implies automatic''
direction of Theorem~\ref{thm:main}, we must explain why the field
of $p$-quasi-automatic series is closed under extraction of roots of
polynomials. The argument we give below implicitly performs a positive
characteristic variant 
of Newton's algorithm; most of the work has already been 
carried out in Chapter~\ref{sec:valued}.
(By contrast, a direct adaptation of
Newton's algorithm to generalized power series gives a transfinite
process; see \cite[Proposition~1]{me2}.)

Before proceeding further, we explicitly check that the class of
$p$-quasi-automatic series is closed under the formation of
Artin-Schreier extensions.

\begin{lemma} \label{L:artin-schreier1}
For any $p$-quasi-automatic $x = \sum x_i t^i \in \FF_q((t^\QQ))$ supported
within $(-\infty, 0)$, there exists
a $p$-quasi-automatic series $y \in \FF_q((t^\QQ))$ such that
$y^{p} - y = x$.
\end{lemma}
\begin{proof}
We can take
\[
y = x^{1/p} + x^{1/p^{2}} + \cdots
\]
once we show that this generalized power series is $p$-quasi-automatic.
There is no loss of generality in assuming that $x t^a$ is $p$-automatic
for some nonnegative integer $a$; in fact, by decimating, we may reduce
to the case where $a=1$.

Write $y = \sum y_i t^i$.
Since 
\[
y_i = x_{ip} + x_{ip^{2}} + \cdots,
\]
for any fixed $j$, the series 
\[
\sum_{i < -p^{-j}} y_i t^i
\]
is $p$-automatic. 
Also, for $i \in [-1, 0)$, the sequence
\[
x_i, x_{i/p}, x_{i/p^{2}}, \dots
\]
is generated by inputting the base $p$ expansions of
$1+i, 1+i/p, 1+i/p^{2}, \dots$ into a fixed finite
automaton; hence there exist integers $m$ and $n$ such that for any
$i \in [-1,0)$, $y_{ip^{-m}} = y_{ip^{-n}}$.

From this we can construct a finite automaton that, upon receiving
$1+i$ at input, returns $y_i$. Namely, the automaticity of 
$\sum_{i < -p^{-n}} y_i t^i$ gives an automaton that returns $y_i$
if $i < -p^{-n}$. If $i \geq -p^{-n}$, then $1+i \geq 1 - p^{-n}$,
so the base $p$ expansion begins with $n$ digits equal to $p-1$,
and conversely. We can thus loop back to wherever we were after
$m$ digits equal to $p-1$; since $y_{ip^{-m}} = y_{ip^{-n}}$,
we end up computing the right coefficient. Hence $y$
is $p$-quasi-automatic, as desired.
\end{proof}

\begin{lemma} \label{L:artin-schreier2}
Let $R_q \subset \FF_q((t^\QQ))$ be the closure, under the valuation $v$,
of the field of $p$-quasi-automatic series. Then for
any $a,b \in R_q$ with $a \neq 0$, the equation 
\[
z^p - az = b
\]
has $p$ distinct solutions in $R_{q'}$ for some power $q'$ of $q$.
\end{lemma}
\begin{proof}
We first check the case $b=0$, i.e., we show that $z^{p-1} = a$ has
$p-1$ distinct solutions in $R_{q'}$ for some $q'$.
Write $a = a_0 t^i (1 + u)$, where $a_0 \in \FF_q$, $i \in \QQ$,
and $v(u) > 0$. Choose $q'$ so that
$a_0$ has a full set of $(p-1)$-st roots in $\FF_{q'}$; we can then
take 
\[
z = a_0^{1/(p-1)} t^{i/(p-1)} \sum_{j=0}^\infty \binom{j}{1/(p-1)} 
u^j
\]
for any $(p-1)$-st root $a_0^{1/(p-1)}$ of $a_0$ .

We now proceed to the case of general $b$. By the above argument,
we can reduce to the case $a=1$. We can then split $b = b_- + b_+$,
where $b_-$ is supported on $(-\infty, 0)$ and
$b_+$ is supported on $[0, \infty)$, and treat the cases
$b = b_-$ and $b = b_+$ separately. The former case
is precisely Lemma~\ref{L:artin-schreier1}. As for the latter case,
let $b_0$ be the constant coefficient of $b_+$, and choose $q'$
so that the equation $z^p - z = b_0$ has distinct
roots $c_1, \dots, c_p \in \FF_{q'}$.
We may then take
\[
z = c_i - \sum_{j=0}^\infty b_0^{p^j}
\]
to obtain the desired solutions.
\end{proof}

\begin{prop} \label{P:Newton}
Let $R_q \subset \FF_q((t^\QQ))$ be the closure, under the valuation $v$,
of the field of $p$-quasi-automatic series, and let
$R$ be the union of the rings $R_{q'}$ over all powers $q'$ of $q$. Then
$R$ (which is also a field) is algebraically closed.
\end{prop}
\begin{proof}
By Lemma~\ref{L:ore}, it suffices to show that for every monic
twisted polynomial
$P(F)$ over $R_q$, the polynomial $P(F)(z)$ factors completely over
$R_{q'}$ for some $q'$. Since $R_q$ is perfect, we may
factor off $F$ on the right if it appears, to reduce to the case where
$P$ has nonzero constant coefficient, or equivalently
(by Lemma~\ref{L:Fp-vector space}) where $P(F)(z)$ has no repeated roots.

By Corollary~\ref{C:splits completely}, we can write
$P = Q_1 \cdots Q_n$ for some monic linear twisted polynomials
$Q_1, \dots, Q_n$ over $R_{q'}$ for some $q'$.
Write $Q_i = F - c_i$ (where $c_i \neq 0$);
the process of finding roots of $P(F)(z)$ can then be described as the
process of finding solutions of the system of equations
\begin{align*}
z_1^p - c_1 z_1 &= 0 \\
z_2^p - c_2 z_2 &= z_1 \\
&\vdots \\
z_n^p - c_n z_n &= z_{n-1}
\end{align*}
(the roots of $P(F)(z)$ being precisely the possible values of $z_n$).
By repeated applications of Lemma~\ref{L:artin-schreier2}, the
system has $p^n$ distinct solutions,
and so $P(F)(z)$ splits completely as desired.
\end{proof}

At this point, we can now deduce that ``algebraic implies
automatic'' in much the same manner as in Proposition~\ref{P:algtoaut1}.
This means in particular that again we will need to invoke Christol's theorem.

\begin{prop} \label{P:algtoaut2}
Let $x = \sum x_i t^i \in \FF_q((t^\QQ))$ be a generalized power
series which is algebraic over $\FF_q(t)$. Then $x$ is
$p$-quasi-automatic.
\end{prop}
\begin{proof}
There exists a 
polynomial $P(z)$ over $\FF_q(t^{1/p^m})$, for some nonnegative
integer $m$, such that $P$ has $x$ as a root with multiplicity 1;
by replacing $x$ with $x^{p^m}$, we may reduce without loss of 
generality to the case where $m = 1$.
Choose $c \in \QQ$ such that $c > v(x-x')$ for any root
$x' \neq x$ of $P$.

By Proposition~\ref{P:Newton}, for any $c \in \QQ$, there exists
a power $q'$ of $q$ and 
a $p$-quasi-automatic series $y$ over $\FF_{q'}((t^\QQ))$
such that $v(x-y) \geq c$. The polynomial $P(z+y)$ then has
exactly one root of slope at least $c$, namely $x-y$.

By Proposition~\ref{P:auttoalg}, $y$ is algebraic over $\FF_q(t)$.
Let $K$ be the finite extension of $\FF_{q'}((t))$ obtained by
adjoining $y$. Then $K$ is complete under $v$; 
by Corollary~\ref{C:slope factorization}, we may split off
the unique factor of $P(z+y)$ of slope at least $c$.
In other words, $x-y \in K$, and so $x \in K$.

At this point the argument parallels that of 
Proposition~\ref{P:algtoaut1}.
Let $m$ be the degree of the minimal polynomial of $y$ over 
$\FF_q(t)$.
For $j=0, \dots, m-1$, write
\[
y^{pj} = \sum_{i=0}^{m-1} a_{ij} y^i
\]
with $a_{ij} \in \FF_q((t))$; then the $a_{ij}$ are algebraic over
$\FF_q(t)$.
Choose $n$ minimal such that $x, x^p, \dots, x^{p^m}$ are linearly dependent
over $\FF_q((t))$.
Write $x^{p^j} = \sum_{i=0}^{m-1} b_{ij} y^i$ with $b_{ij} \in \FF_q((t))$;
we then have the equations
\[
b_{i(j+1)} = \sum_{l=0}^{m-1} b_{lj} a_{il} \qquad (j=0, \dots, n-1).
\]
Let $c_0 x + \cdots + c_n x^{p^n} = 0$ be a linear relation
over $\FF_q((t))$, which we may normalize by setting $c_0 = 1$; then
the $c_i$ are algebraic over $\FF_q(t)$.
Also, writing $x = -c_1 x^p - \cdots - c_n (x^{p^{n-1}})^p$, we have
\begin{align*}
\sum_{i=0}^{m-1} b_{i0} y^i &= \sum_{j=0}^{n-1} -c_{j+1} (x^{p^j})^p \\
&= \sum_{i=0}^{m-1} \sum_{j=0}^{n-1} -c_{j+1} b_{lj}^p y^{pl} \\
&= \sum_{i=0}^{m-1} \sum_{j=0}^{n-1} \sum_{l=0}^{m-1}
-c_{j+1} b_{lj}^p a_{il} y^i,
\end{align*}
and hence
\[
b_{i0} = \sum_{j=0}^{n-1} \sum_{l=0}^{m-1} -c_{j+1} b_{lj}^p a_{il}.
\]
We thus have a system of equations in the $b_{ij}$
as in Lemma~\ref{lem:vector} (with $B$ equal to the identity matrix);
this allows us to conclude that each $b_{ij}$ is algebraic over
$\FF_q(t)$.

By Christol's theorem (Theorem~\ref{T:christol}), each $b_{ij} \in \FF_q((t))$
is $p$-automatic. Hence $x = \sum_{i=0}^{m-1} b_{i0} y^i$
is $p$-quasi-automatic, as desired.
\end{proof}

\section{Further questions}
\label{sec:further}

We end with some further questions about automata and generalized
power series, of both an algorithmic and a theoretical nature.

\subsection{Algorithmics of automatic series}
\label{subsec:algor}

It seems that one should be able to make explicit calculations
in the algebraic closure of $\FF_q(t)$ using automatic power series;
we have not made any systematic attempt to do so, but it is worth
recording some observations here for the benefit of anyone considering
doing so in the future.

The issue of computing in the algebraic closure of $\FF_q(t)$ using
generalized power series 
resembles that of computing in the algebraic closure of $\QQ$
using complex approximations of algebraic numbers, and
some lessons may be profitably drawn from that case. In particular,
it may be worth computing ``approximately'' and not exactly with
automatic power series, in a form of ``interval arithmetic''. 

In contrast
with the $\QQ$-analogue, however, there are two ways to truncate a
computation with automata: one can ignore large powers of $t$ (analogous
to working with a complex approximation of an algebraic number), but one
can also prune the automata by ignoring states that cannot be reached
in some particular number of steps from the initial state.
The latter may be crucial for making multiplication of automatic series
efficient, as the methods we have described seem to entail exponential
growth in the number of states over the course of a sequence of
arithmetic operations. Further analysis will be needed, however, to
determine how much pruning one can get away with, and how easily
one can recover the missing precision in case it is needed again.

Our techniques seem to depend rather badly on the size of the finite
field $\FF_q$ under consideration, but it is possible this dependence
can be ameliorated. For instance, some of this dependence (like the
complexity of a product) really depend not on $q$ but on the characteristic
$p$, and so are not so much of an issue when working over a field of
small characteristic (as long as one decomposes everything over a basis
of $\FF_q$ over $\FF_p$). Also, it may be possible to work even in
large characteristic by writing everything in terms of a small
additive basis of the finite field (e.g., for $\FF_p$, use the powers
of 2 less than $p$), at least if one is willing to truncate as in the
previous paragraph.

One additional concern that arises when the characteristic is not
small is that our method for extracting roots of polynomials requires
working with additive polynomials. Given an ordinary polynomial $P(z)$ 
of degree
$n$, one easily obtains an additive polynomial of degree at most $p^n$
which has $P(z)$ as a factor (by reducing $z, z^p, z^{p^2},
\dots$ modulo $P$); however, the complexity of the coefficients of the
new polynomial grows exponentially in $n$. It would be of some interest
to develop a form of Newton's algorithm to deal directly with
ordinary polynomials.

\subsection{Multivariate series and automata}
\label{subsec:multi}

We conclude by mentioning a multivariate version of Christol's theorem
and conjecturing a generalized power series analogue.
Following \cite[Chapter~14]{as}, we restrict
our notion of ``multivariate'' to ``bivariate'' for notational
simplicity, and leave it to the reader's imagination to come up with
full multivariate analogues.

For $b$ an integer greater than 1,
a \emph{valid pair of base $b$ expansions}
is a pair of strings $(s_1 \dots s_n, t_1 \dots t_n)$ of equal length
over the alphabet $\Sigma = \{0, \dots, b-1, .\}$ such that
$s_1$ and $t_1$ are not both 0, $s_n$ and $t_n$ are not both 0,
each of the strings $s_1\dots s_n$ and $t_1 \dots t_n$ contains exactly
one radix point, and those radix points occur at the same index $k$.
We define the \emph{value} of such a pair to be the pair
\[
\left(\sum_{i=1}^{k-1} s_i b^{k-1-i} + \sum_{i=k+1}^n s_i b^{k-i},
\quad \sum_{i=1}^{k-1} t_i b^{k-1-i} + \sum_{i=k+1}^n t_i b^{k-i}\right);
\]
then the value function gives a bijection between the set of valid
pairs of base $b$ expansions and $S_b \times S_b$. Let $s$ denote the
inverse function.
(Note that it may happen that one string or the other has some leading
and/or trailing zeroes, since we are forcing them to have the same length.)

We may identify pairs of strings over $\Sigma_b$ of equal length with strings
over $\Sigma_b \times \Sigma_b$ in the obvious fashion, and under this
identification, the set of valid pairs of base $b$ expansions are seen to
form a regular language. With that in mind,
we declare
a function $f: S_b \times S_b \to \Delta$ to be \emph{$b$-automatic}
if there exists a DFAO with input alphabet $\Sigma = \Sigma_b \times \Sigma_b$
and output alphabet $\Delta$ such that for any pair $(v, w) \in S_b
\times S_b$, $f(v,w) = f_M(s(v,w))$. We declare a double sequence
$\{a_{i,j}\}_{i,j=0}^\infty$ over $\Delta$ to be $b$-automatic if for 
some $\star \notin \Delta$, the function
$f: S_b \times S_b \to \Delta \cup \{\star\}$ given by
\[
f(v,w) = \begin{cases} a_{v,w} & v,w \in \ZZ \\
0 & \mbox{otherwise}
\end{cases}
\]
is $b$-automatic.

In this language, the bivariate version of Christol's theorem is the
following result, due to
Salon \cite{salon1}, \cite{salon2} (see also 
\cite[Theorem~14.4.1]{as}). Here $\FF_q(t,u)$ denotes the
fraction field of the polynomial ring $\FF_q[t,u]$; this field is contained
in the fraction field of the power series ring $\FF_q \llbracket t,u
\rrbracket$.
\begin{theorem} \label{thm:twodim}
Let $q$ be a power of the prime number $p$, and let
$\{a_{i,j}\}_{i,j=0}^\infty$ be a double sequence over $\FF_q$. Then
the double series $\sum_{i,j=0}^\infty a_{i,j} t^i u^j \in 
\FF_q \llbracket t,u \rrbracket$ is algebraic over $\FF_q(t,u)$
if and only if the double sequence $\{a_{i,j}\}_{i,j=0}^\infty$
is $p$-automatic.
\end{theorem}

To even formulate a generalized power series analogue 
of Theorem~\ref{thm:twodim}, we must decide what we mean by ``generalized
power series in two variables''. The construction below is natural enough,
but we are not aware of a prior appearance in the literature.

Let $G$ be a \emph{partially ordered} abelian group (written additively) with
identity element $0$; that is, $G$ is an abelian group equipped with
a binary relation $>$ such that for all $a,b,c \in G$,
\begin{gather*}
a \not> a \\
a > b, b > c \Rightarrow a > c \\
a > b \Leftrightarrow a+c > b+c.
\end{gather*}
Let $P$ be the set of $a \in G$ with $a > 0$; $P$ is again called the
\emph{positive cone} of $G$. We write $a \geq b$ to mean $a > b$ or $a = b$,
and $a \leq b$ to mean $b \geq a$. Then one has the following analogue
of Lemma~\ref{lem:wo}, whose proof we leave to the reader.
\begin{lemma}\label{lem:wpo}
Let $S$ be a subset of $G$. Then the
following two conditions are equivalent.
\begin{enumerate}
\item[(a)] Every nonempty subset of $S$ has at least one, but only
finitely many, minimal elements. (An element $x \in S$ is \emph{minimal}
if $y \in S$ and $y \leq x$ imply $y = x$. Note that ``minimal'' does not
mean ``smallest''.)
\item[(b)] Any sequence $s_1, s_2, \dots$ over $S$ contains an infinite
weakly increasing subsequence $s_{i_1} \leq s_{i_2} \leq \cdots$.
\end{enumerate}
\end{lemma}
A subset $S$ of $G$ is \emph{well-partially-ordered} 
if it satisfies either of the
equivalent conditions of Lemma~\ref{lem:wpo}. 
(This condition has cropped up repeatedly in combinatorial situations;
see \cite{kruskal} for a survey, if a somewhat dated one.)
Then for any ring $R$,
the set of functions $f: G \to R$ which have well-partially-ordered
support forms a ring under termwise addition and convolution.

For $G = \QQ \times \QQ$, let $R((t^\QQ, u^\QQ))$ denote the ring just
constructed; we refer to its elements as ``generalized double
Laurent series''. Then by analogy with Theorem~\ref{thm:main}, we formulate
the following conjecture.
\begin{conj} \label{conj:twodim}
Let $q$ be a power of the prime $p$, and
let $f: \QQ \times \QQ \to \FF_q$ be a function whose support $S$
is well-quasi-ordered.
Then the corresponding generalized double Laurent series $\sum_{i,j} f(i,j) 
t^i u^j
\in \FF_q (( t^\QQ, u^\QQ ))$ 
is algebraic over $\FF_q(t,u)$ if and only if the
following conditions hold. 
\begin{enumerate}
\item[(a)]
For some positive integers $a$ and $b$, the set
$aS+b = \{(ai+b, aj+b): (i,j) \in S\}$ is contained in $S_p \times S_p$.
\item[(b)]
For some $a,b$ for which (a) holds, the function
$f_{a,b}: S_p \times S_p \to \FF_q$ given by 
$f_{a,b}(x,y) = f((x-b)/a, (y-b)/a)$ is $p$-automatic.
\end{enumerate}
Moreover, if these conditions hold, then $f_{a,b}$ is $p$-automatic
for any nonnegative integers $a,b$ for which (a) holds.
\end{conj}

An affirmative answer to Conjecture~\ref{conj:twodim} would imply that
if the generalized double Laurent series $\sum_{i,j} c_{i,j} t^i u^j$ is
algebraic over $\FF_q(t,u)$, then the diagonal series
$\sum_i c_{i,i} t^i$ is algebraic over $\FF_q(t)$; for ordinary Laurent
series, this result is due to Deligne \cite{deligne}.

Deligne's result actually allows an arbitrary field of positive characteristic
in place of $\FF_q$, but his proof is in the context of sophisticated
algebro-geometric machinery (vanishing cycles in \'etale cohomology); 
a proof of Deligne's general result
in the spirit of automata-theoretic methods was given
by Sharif and Woodcock \cite{sw}. It may be possible to state and prove
an analogous assertion in the generalized Laurent series setting by
giving a suitable multivariate extension of the results of \cite{me}.

\end{document}